# Estimators of different *delta* coefficients based on the unbiased estimator of the expected proportions of agreements

*Martín Andrés, A.[1]\* and Álvarez Hernández, M.[2,3]*

[1] Biostatistics. Emeritus Professor of the University of Granada, Granada, Spain.
[2] Departamento de Estadística e Investigación Operativa, Universidad de Vigo, Spain.
[3] CITMAga, 15782 Santiago de Compostela, La Coruña, Spain.

## ABSTRACT

To measure the degree of agreement between two observers that independently classify $n$ subjects within $K$ categories, it is common to use different *kappa* type coefficients, the most common of which is the $\kappa_C$ coefficient (Cohen's *kappa*). As $\kappa_C$ has some weaknesses -such as its poor performance with highly unbalanced marginal distributions-, the $\Delta$ coefficient is sometimes used, based on the *delta* response model. This model allows us to obtain other parameters like: (a) the $\alpha_i$ contribution of each $i$ category to the value of the global agreement $\Delta=\Sigma\alpha_i$; and (b) the *consistency* $\mathcal{S}_i$ in the category $i$ (degree of agreement in the category $i$), a more appropriate parameter than the *kappa* value obtained by collapsing the data into the category $i$. It has recently been shown that the classic estimator $\hat{\kappa}_C$ underestimates $\kappa_C$, having obtained a new estimator $\hat{\kappa}_{CU}$ which is less biased. This article demonstrates that something similar happens to the known estimators $\hat{\Delta}$, $\hat{\alpha}_i$, and $\hat{\mathcal{S}}_i$ of $\Delta$, $\alpha_i$, and $\mathcal{S}_i$ (respectively), proposes new and less biased estimators $\hat{\Delta}_U$, $\hat{\alpha}_{iU}$, and $\hat{\mathcal{S}}_{iU}$, determines their variances, analyses the behaviour of all estimators, and concludes that the new estimators should be used when $n$ or $K$ are small (at least when $n\leq 50$ or $K\leq 3$). Additionally, the case where one of the raters is a gold standard is contemplated, in which situation two new parameters arise: the *conformity* (the rater's capability to recognize a subject in the category $i$) and the *predictivity* (the reliability of a response $i$ by the rater).

\* Email: amartina@ugr.es. Phone: 34-58-244080.

# 1. Introduction

Let there be two raters that independently classify $n$ subjects in $K$ nominal categories. Given a subject, Observer 1 classifies it as type $i$ ($i = 1, 2, ..., K$) and Observer 2 as type $j$ ($j = 1, 2, ..., K$). This leads to a table of absolute frequencies $x_{ij}$ like that of Table 1(a), in which $x_{i\bullet}=\Sigma_j x_{ij}$ and $x_{\bullet j}=\Sigma_i x_{ij}$. The random variable $\{x_{ij}\}$ follows a multinomial distribution $M\{n \mid p_{ij}\}$, where $p_{ij}$ is the probability of a subject being classified in the cell $(i, j)$: see Table 1(b), in which $p_{i\bullet}=\Sigma_j p_{ij}$ and $p_{\bullet j}=\Sigma_i p_{ij}$ are the marginal distributions. In this scenario, the primary interest of researchers is to determine the degree of concordance or agreement between both raters. From a population point of view, this degree of agreement cannot be the populational index of observed agreements $I_o=\Sigma_i p_{ii}$ -the sum of all the probabilities of agreement -, since some of these agreements may be random. This is why it is standard procedure to eliminate the random effect defining a *kappa* type coefficient $\kappa = (I_o - I_e)/(1 - I_e)$, where $I_e$ is the populational index of expected agreements -the sum of the probabilities of agreements that happen randomly - and $\kappa$ is the degree of population agreement. Depending on how $I_e$ is defined we obtain different *kappa* -$\kappa_S$ coefficients (Scott 1955), $\kappa_C$ (Cohen 1960), $\kappa_K$ (Krippendorff 2004), and $\kappa_G$ (Gwet 2008), among others-, all of which are based on the marginal distributions. The standard procedure is to define $I_e=\Sigma_i p_{i\bullet} p_{\bullet i}$, thus obtaining the $\kappa_C$ value of Cohen.

The estimation of these $\kappa$ coefficients has the general form of $\hat{\kappa} = \left(\hat{I}_o - \hat{I}_e\right)\big/\left(1 - \hat{I}_e\right)$, where the values $\hat{\kappa}$, $\hat{I}_o = \Sigma \bar{p}_{ii}$ (with $\bar{p}_{ij}=x_{ij}/n$) and $\hat{I}_e$ are the sample estimations of the previous population parameters $\kappa$, $I_o$, and $I_e$, respectively. In the case of $\hat{\kappa}_C$ it is $\hat{I}_e = \Sigma_i \bar{p}_{i\bullet} \bar{p}_{\bullet i}$, where $\bar{p}_{i\bullet}=\Sigma_j \bar{p}_{ij}$ and $\bar{p}_{\bullet j} = \Sigma_i \bar{p}_{ij}$. Something similar happens with the other coefficients $\hat{\kappa}_S$, $\hat{\kappa}_K$, and $\hat{\kappa}_G$. Nevertheless, all of the estimators $\hat{I}_e$ -of any of the definitions of $I_e$ which were proposed – are based on the assumption that $\bar{p}_{i\bullet} \bar{p}_{\bullet i}$ is an unbiased estimator of $p_{i\bullet} p_{\bullet i}$. Martín Andrés and Álvarez Hernández (2025) recently observed that this is not true -since the expected value of a product is not generally the product of the expected values -, and they proposed the necessary expressions to correct that bias, thus improving the estimation of $\kappa$.

Although $\hat{\kappa}_C$ is a very popular measure of agreement, it has two very significant limitations (Brennan and Prediger 1981, Agresti *et al.* 1995, Guggenmoos-Holzmann and Vonk 1998, Nelson and Pepe 2000, Martín Andrés and Femia Marzo 2004 and 2005, and Erdmann *et al.* 2015): its dependence on marginal distributions and its difficulty in evaluating the degree of agreement in each category. Martín Andrés and Femia Marzo (2004 and 2005) and Martín Andrés and Álvarez Hernández (2022) proposed a response model (*delta* model) that gives rise to the agreement measure *delta* ($\Delta$) which does not have the limitations of $\kappa_C$ (Ato *et al.* 2011, Shankar and Bangdiwala 2014, Giammarino *et al.* 2021). For instance, with $K=2$ and data $x_{11}=80$, $x_{12}=10$, $x_{21}=10$ y $x_{22}=0$ by Nelson and Pepe (2000), $\hat{\kappa}_C=-0.11$ and $\hat{\Delta}=0.58$ are obtained. In the classic definitions of agreement, the $p_{ij}$ probabilities are not submitted to any model, but in the case of the *delta* model (with the notation of Martín Andrés and Álvarez Hernández 2022)

$$p_{ij} = \delta_{ij}\alpha_i + B\pi_{i1}\pi_{j2} \text{ with } B = 1-\Delta, \ \Delta=\Sigma_i\alpha_i, \ 0\leq\pi_{ir}\leq1 \text{ and } \Sigma_i\pi_{ir}=1 \ (\forall \ r=1, 2), \quad \textbf{(1)}$$

when $\delta_{ij}$ are the Kronecker deltas. In this model, $\Delta$ is the total proportion of non-random agreements (global agreement), $\alpha_i$ is the proportion of agreements in Category $i$ which are not random (contribution of category $i$ to the global agreement), and $\{\pi_{ir}\}$ is the distribution of the response that the rater $r$ makes randomly (in the proportion of subjects $B=1-\Delta$ that classifies by chance). Note that the model indicates that $p_{ii}=\alpha_i+B\pi_{i1}\pi_{i2}$ for $i=j$, so the proportion $p_{ii}$ of observed agreements (theoretical) in the category $i$ decomposes into two addends: the proportions of agreements NOT due by chance ($\alpha_{ii}$) y those YES due by chance ($B\pi_{i1}\pi_{i2}$). That is why the global raw agreement is $\Sigma_i p_{ii}$, while the global agreement corrected by chance is

$\Delta = \Sigma_i \alpha_i$.

The model has the advantage of providing results similar to those of the classic $\kappa_C$ in normal circumstances, but provides a greater degree of agreement than $\kappa_C$ when the marginal distributions are very unbalanced in the same direction. Coefficient $\Delta$ can also be put in the format of $\kappa$ since, due to expression (1), $\Delta = (I_o - I_\pi)/(1 - I_\pi)$ with $I_\pi = \Sigma_i \pi_{i1} \pi_{i2}$. The aforementioned authors provide the estimation of the parameters of the *delta* model through the method of maximum likelihood -values $\hat{\alpha}_i = \overline{p}_{ii} - \left(1 - \hat{\Delta}\right) \hat{\pi}_{i1} \hat{\pi}_{i2}$, $\hat{\pi}_{ir}$ and $\hat{\Delta} = \Sigma_i \hat{\alpha}_i$ (see Appendix A)-, and it also occurs that $\hat{\Delta} = \left(\hat{I}_o - \hat{I}_\pi\right) \big/ \left(1 - \hat{I}_\pi\right)$, where $\hat{I}_\pi = \Sigma_i \hat{\pi}_{i1} \hat{\pi}_{i2}$. As once again $\hat{\pi}_{i1} \hat{\pi}_{i2}$ is not an unbiased estimator of $\pi_{i1} \pi_{i2}$, then our first objective will be to correct this defect and improve estimator $\hat{\Delta}$.

Additionally, researchers can also have two secondary interests. Firstly, to determine the contribution $\alpha_i$ of category $i$ to the global agreement $\Delta$. Secondly, to determine the degree of agreement in category $i$ which is not random i.e. the *consistency* $\mathcal{S}_i$ in category $i$

$$\mathcal{S}_i = \frac{2\alpha_i}{p_{i\bullet} + p_{\bullet i}} = \frac{2\alpha_i}{2\alpha_i + \left(1 - \Delta\right)\left(\pi_{i1} + \pi_{i2}\right)}. \qquad (2)$$

This parameter, which comes from the raw agreement index in category $i$ by Cicchetti and Feinstein (1990) -$S_i = 2p_{ii}/(p_{i\bullet} + p_{\bullet i})$, in which $p_{ii}$ has been replaced by $\alpha_i$- refers to the proportion of (non-random) agreements between all of the subjects classified in category $i$ by either of the two raters, and was analysed and estimated in the context of the *delta* model by Martín Andrés and Femia Marzo (2005). These authors demonstrated that $\mathcal{S}_i$ has the same objective as that sought when defining $\kappa_C$ for the collapsed data in the category $i$ -parameter $\kappa_{C(i)}$-, an objective which is always achieved with $\mathcal{S}_i$, but which sometimes is not achieved with $\kappa_{C(i)}$. In fact, a parameter of global agreement based on the collapsed data in the category $i$ measures the degree of global agreement in the new situation, but does not measure the degree of agreement in category *i* since, for the same reason, it should also measure the degree of agreement in the category "no *i*". As in both cases ($\alpha_i$ and $\mathcal{S}_i$) the estimation of the parameter ($\hat{\alpha}_i$ and $\hat{\mathcal{S}}_i$) depends on the estimation of the product $\pi_{i1}\pi_{i2}$, the conclusion is the same as in the previous paragraph: the two estimations are biased and that bias should be corrected as much as possible. This will be the second objective of this article.

The third objective is obvious: to determine the variance of each of the new estimators obtained.

The fourth and last objective, discussed in a more superficial way, consists of measuring the quality of a rater's response when the other rater is a gold standard, which is related to diagnostic methods.

## 2. Approximately unbiased estimator of $\pi_{i1}\pi_{i2}$ and new estimators proposed for $\Delta$, $\alpha_i$ and $S_i$ when $K > 2$.

In the following it is assumed that $\hat{\alpha}_i$, $\hat{\pi}_{ir}$, and $\hat{\Delta}$ are the maximum likelihood estimators of $\alpha_i$, $\pi_{ir}$, and $\Delta$ respectively, obtained under to the *delta* model. Appendix A summarizes the procedure to obtain these estimators, and at http://www.ugr.es/local/bioest/ software (Section "Agreement Among Raters") it is possible to download two free programmes (Delta.exe and Multi_Rater_Delta.exe) to run this model. Appendix A also demonstrates that

$$E\left(\hat{\pi}_{i1}\hat{\pi}_{i2}\right) = \pi_{i1}\pi_{i2} + E_i \text{ where } E_i = \frac{1}{n\left(1-\Delta\right)}\left\{\pi_{i1}\pi_{i2} - \frac{X_i\left(X-X_i\right)}{X-1}\right\}, \qquad \textbf{(3)}$$

with $X = \Sigma_i X_i$ and $X_i = \pi_{i1}\pi_{i2}/(\pi_{i1}+\pi_{i2}-1)$. This means that $\hat{\pi}_{i1}\hat{\pi}_{i2}$ is a biased estimator of $\pi_{i1}\pi_{i2}$ and that it is not possible to easily correct the bias of $\hat{\pi}_{i1}\hat{\pi}_{i2}$. The best option (see Appendix A) is the following approximately unbiased estimation:

$$\left(\widehat{\pi_{i1}\pi_{i2}}\right)_U = \hat{\pi}_{i1}\hat{\pi}_{i2} - \hat{E}_i \text{ where } \hat{E}_i = \frac{1}{n\left(1-\hat{\Delta}\right)}\left\{\hat{\pi}_{i1}\hat{\pi}_{i2} - \frac{\hat{X}_i\left(\hat{X}-\hat{X}_i\right)}{\hat{X}-1}\right\}, \qquad \textbf{(4)}$$

where $\hat{X} = \sum_i \hat{X}_i$ and $\hat{X}_i = \hat{\pi}_{i1}\hat{\pi}_{i2} \big/ \left(\hat{\pi}_{i1} + \hat{\pi}_{i2} - 1\right)$.

Expression (4) also causes the traditional estimators of $I_\pi$, $\Delta$, $\alpha_i$, and $\mathcal{S}_i$ to change, since all of them depend on this expression. As the classic estimators are:

$$\hat{I}_\pi = \Sigma \, \hat{\pi}_{i1}\hat{\pi}_{i2}, \quad \hat{\Delta} = \frac{\hat{I}_0 - \hat{I}_\pi}{1 - \hat{I}_\pi}, \quad \hat{\alpha}_i = \overline{p}_{ii} - \left(1 - \hat{\Delta}\right)\hat{\pi}_{i1}\hat{\pi}_{i2} \text{ and } \hat{\mathcal{S}}_i = \frac{2\hat{\alpha}_i}{\overline{p}_{i\bullet} + \overline{p}_{\bullet i}}, \tag{5}$$

the new (approximately unbiased) estimators will be,

$$\hat{I}_{\pi U} = \hat{I}_\pi - \hat{E}, \text{ where } \hat{E} = \sum_i \hat{E}_i = \frac{1}{n\left(1 - \hat{\Delta}\right)}\left\{\hat{I}_\pi - \frac{\hat{X}^2 - \sum_i \hat{X}_i^2}{\hat{X} - 1}\right\}, \tag{6}$$

$$\hat{\Delta}_U = \frac{\hat{I}_0 - \hat{I}_{\pi U}}{1 - \hat{I}_{\pi U}}, \quad \hat{\alpha}_{iU} = \overline{p}_{ii} - \left(1 - \hat{\Delta}_U\right)\left(\hat{\pi}_{i1}\hat{\pi}_{i2} - \hat{E}_i\right) \text{ and } \hat{\mathcal{S}}_{iU} = \frac{2\hat{\alpha}_{iU}}{\overline{p}_{i\bullet} + \overline{p}_{\bullet i}}, \tag{7}$$

where it can be observed that the definitions are coherent: $\Sigma_i \, \hat{\alpha}_{iU} = \hat{\Delta}_U$. In general it can be expected that most of the time $\hat{\Delta}_U \geq \hat{\Delta}$, $\hat{\alpha}_{iU} \geq \hat{\alpha}_i$, and $\hat{\mathcal{S}}_{iU} \geq \hat{\mathcal{S}}_i$. The reason for this is that it will generally happen that $\hat{\pi}_{i1} + \hat{\pi}_{i2} - 1 < 0$, so that $\hat{X}_i < 0$ and, through expressions (4) and (6), $\hat{E}_i > 0$, $\left(\widehat{\pi_{i1}\pi_{i2}}\right)_U < \hat{\pi}_{i1}\hat{\pi}_{i2}$, and $\hat{I}_{\pi U} \leq \hat{I}_\pi$. The simulations in Section 4 confirm this question.

## 3. Variances of the new estimators.

Martín Andrés and Femia Marzo (2004 and 2005) provided the following asymptotic variances of the estimators $\hat{\Delta}$, $\hat{\alpha}_i$ and $\hat{\mathcal{S}}_i$, which we express here -and in everything that follows- in the format of Martín Andrés and Álvarez Hernández (2022),

$$V_A\left(\hat{\Delta}\right) = \frac{1 - \Delta}{n}\left\{\Delta + \frac{X}{X - 1}\right\}, \quad V_A\left(\hat{\alpha}_i\right) = \frac{H_i + \alpha_i\left(1 - \alpha_i\right)}{n}, \text{ and}$$

$$V_A\left(\hat{\mathcal{S}}_i\right) = \frac{4H_i + \mathcal{S}_i\left\{2t_i - 3t_i\mathcal{S}_i + 2p_{ii}\mathcal{S}_i\right\}}{nt_i^2}, \tag{8}$$

where $t_i = p_{i\bullet} + p_{\bullet i}$ and

$$H_i = \left(1 - \Delta\right)X_i\left\{\frac{X_i}{X - 1} - 1\right\}, \tag{9}$$

as well as the following estimated variances:

$$\hat{V}\left(\hat{\Delta}\right)=\frac{1-\hat{\Delta}}{n}\left\{\hat{\Delta}+\frac{\hat{X}}{\hat{X}-1}\right\},\ \hat{V}\left(\hat{\alpha}_i\right)=\frac{\hat{H}_i+\hat{\alpha}_i\left(1-\hat{\alpha}_i\right)}{n},\ \text{and}$$

$$\hat{V}\left(\hat{\mathcal{S}}_i\right)=\frac{4\hat{H}_i+\hat{\mathcal{S}}_i\left\{2\overline{t_i}-3\overline{t_i}\hat{\mathcal{S}}_i+2\,\overline{p}_{ii}\hat{\mathcal{S}}_i\right\}}{n\overline{t_i}^2},\tag{10}$$

where $\overline{t_i}=\overline{p}_{i\bullet}+\overline{p}_{\bullet i}$ and

$$\hat{H}_i=\left(1-\hat{\Delta}\right)\hat{X}_i\left\{\frac{\hat{X}_i}{\hat{X}-1}-1\right\}.\tag{11}$$

The determination of the formula of the variance of the new estimators $\hat{\Delta}_U$, $\hat{\alpha}_{iU}$, and $\hat{\mathcal{S}}_{iU}$ is complicated and its expression will be very complex. Given that the extra terms of these estimators - $\hat{E}_i$ and $\hat{E}$ of expressions (4) and (6) respectively- are divided by $n$, it is reasonable to assume that the formulas will be similar to those of expressions (10), with the appropriate changes. For example, as $(1-\hat{\Delta}) = (1-\hat{\Delta}_U)\{1+\hat{E}/(1-\hat{I}_\pi)\} = (1-\hat{\Delta}_U)+\hat{B}$, where $\hat{B}=(1-\hat{\Delta}_U)\hat{E}/(1-\hat{I}_\pi)$, then $V(1-\hat{\Delta}) = V(1-\hat{\Delta}_U)+V(\hat{B})-2\text{Cov}\{\hat{\Delta}_U,\hat{B}\}\approx V(1-\hat{\Delta}_U)$ -since $\hat{B}$ is $O(n^{-1})$- and therefore $V(\hat{\Delta})\approx V(\hat{\Delta}_U)$. That is why the proposed variances are given by that the following modifications of expressions (10) and (11):

$$\hat{V}\left(\hat{\Delta}_U\right)=\frac{1-\hat{\Delta}_U}{n}\left\{\hat{\Delta}_U+\frac{\hat{X}}{\hat{X}-1}\right\},\ \hat{V}\left(\hat{\alpha}_{iU}\right)=\frac{\hat{H}_{iU}+\hat{\alpha}_{iU}\left(1-\hat{\alpha}_{iU}\right)}{n}\ \text{y}$$

$$\hat{V}\left(\hat{\mathcal{S}}_{iU}\right)=\frac{4\hat{H}_{iU}+\hat{\mathcal{S}}_{iU}\left\{2\overline{t_i}-3\overline{t_i}\hat{\mathcal{S}}_{iU}+2\,\overline{p}_{ii}\hat{\mathcal{S}}_{iU}\right\}}{n\overline{t_i}^2},\tag{12}$$

where:

$$\hat{H}_{iU}=\left(1-\hat{\Delta}_U\right)\hat{X}_i\left\{\frac{\hat{X}_i}{\hat{X}-1}-1\right\}.\tag{13}$$

**4. Simulations.**

This Section will determine the bias of the three classic estimators ($\hat{\Delta}$, $\hat{\alpha}_i$, and $\hat{\mathcal{S}}_i$) and of the three new estimators ($\hat{\Delta}_U$, $\hat{\alpha}_{iU}$, and $\hat{\mathcal{S}}_{iU}$), assessing and comparing the behaviour of all of them. Furthermore, this section will study the behaviour of the different formulas of the variances indicated in the previous section. For this purpose, we will consider the 48 setting indicated in Table 2, each one of which refer to a multinomial distribution M$\{p_{ij}; n\}$ whose parameters under the *delta* model are those indicated on each line: values of $n$, $K$, $\alpha_i$, $\pi_{i1}$, and $\pi_{i2}$.

*4.1. The case of parameter $\Delta$.*

To assess the estimators of $\Delta$ and of their variances, the procedure is as follows:

(1) For each setting in Table 2, calculate $V_A\left(\hat{\Delta}\right)$ of the first expression (8) and the value of $\Delta$.

(2) Obtain $N$=10,000 random samples and, in each sample $h$ -with $h = 1, 2, \ldots, N$-, obtain the estimations $\hat{\Delta}_h$ and $\hat{\Delta}_{Uh}$ of expressions (5) and (7), and the values $\hat{V}_h\left(\hat{\Delta}\right)$ and $\hat{V}_h\left(\hat{\Delta}_U\right)$ of expressions (10) and (12).

(3) Calculate the average and the sample variance (denominator $N-1$) of the $N$ values of $\hat{\Delta}_h$ and $\hat{\Delta}_{Uh}$: $\left\{\overline{\hat{\Delta}}, V_E\left(\hat{\Delta}\right)\right\}$ and $\left\{\overline{\hat{\Delta}_U}, V_E\left(\hat{\Delta}_U\right)\right\}$ respectively.

(4) Calculate the average of the $N$ values $\hat{V}_h\left(\hat{\Delta}\right)$ and $\hat{V}_h\left(\hat{\Delta}_U\right)$: $\overline{\hat{V}\left(\hat{\Delta}\right)}$ and $\overline{\hat{V}\left(\hat{\Delta}_U\right)}$ respectively.

(5) With the data obtained, Table 3 is constructed.

For the assessment of the estimators of $\Delta$, it is necessary to take into account that the comparison of $\overline{\hat{\Delta}}$ and $\overline{\hat{\Delta}_U}$ with the known value $\Delta$ allows us to assess the bias of each estimator. It is observed that the two estimators always underestimate the true value $\Delta$, except for the estimator $\hat{\Delta}_U$ which slightly overestimates it on one occasion and coincides with it on another. The average underestimation is 0.044 for $\hat{\Delta}$ and 0.021 for $\hat{\Delta}_U$. Moreover, $\overline{\hat{\Delta}_U} \geq \overline{\hat{\Delta}}$, with a minimum of 0.003 and a maximum of 0.070. It can therefore be seen that $\hat{\Delta}_U$ is clearly a better estimator of $\Delta$ than $\hat{\Delta}$. In a more detailed way, the average of the differences $\hat{\Delta}_U - \hat{\Delta}$: (a) decreases with $n$, and is 0.04, 0.02, and 0.01 for $n = 30$, 50 and 100 respectively, which means that it is necessary to make the correction at least for $n \leq 50$; (b) it decreases with $K$, and is 0.03 and 0.01 for $K=3$ and 5 respectively, which means that it is necessary to make the correction at least for $K \leq 3$; (c) it increases very slightly with $\Delta$, passing from 0.021 in $\Delta=0.4$ to 0.025 in $\Delta=0.8$, which means that that the correction is necessary whether $\Delta$ is small or large. The consequence is that $\hat{\Delta}_U$ is a better estimator than $\hat{\Delta}$ and it should be used when $n$ or $K$ are small. For example, one of the extreme cases is that of the first setting ($n=30$, $K=3$ and $\Delta=0.40$) in which $\overline{\hat{\Delta}}=0.31$ and $\overline{\hat{\Delta}_U}=0.38$. If the maximum differences $\overline{\hat{\Delta}_U} - \overline{\hat{\Delta}}$ are determined in all of the settings with $(n, K) = (3, 30)$, $(3, 50)$, $(3, 100)$, $(5, 30)$, $(5, 50)$, and $(5, 100)$, it is observed that these differences can reach the values 0.08, 0.06, 0.03, 0.06, 0.02, and 0.00, respectively. If it is considered that a difference $\overline{\hat{\Delta}_U} - \overline{\hat{\Delta}}$ of 0.02 is relevant, then we should use $\hat{\Delta}_U$ if $n \leq 50$ or $K \leq 3$. Although this has only been demonstrated with very limited values of $n$ and $K$, and in some specific settings, in practice the demands may be rather less strict, since the values $\overline{\hat{\Delta}_U} - \overline{\hat{\Delta}}$ are average differences and, therefore, the values of interest $\hat{\Delta}_U - \hat{\Delta}$ may be much larger.

For the assessment of the variances, it is necessary to take into account that:

(A) The two estimators $\hat{\Delta}$ and $\hat{\Delta}_U$ have exact unknown variances, but they can be estimated in quite a precise manner through the sample variances $V_E\left(\hat{\Delta}\right)$ and $V_E\left(\hat{\Delta}_U\right)$ described previously in paragraph (3). These values will be referred to from now on as "exact variances", serving as a reference to assess the average of the estimated values: $\overline{\hat{V}\left(\hat{\Delta}\right)}$ and $\overline{\hat{V}\left(\hat{\Delta}_U\right)}$.

(B) The true value of the asymptotic variance $\hat{\Delta}$ is given by the value $V_A\left(\hat{\Delta}\right)$ of expression (8), which is not an estimated value, since it is obtained in the true values of the parameters in each setting. The comparison of $V_A\left(\hat{\Delta}\right)$ (exact asymptotic variance) with $V_E\left(\hat{\Delta}\right)$ (exact variance) will indicate how good the asymptotic formula actually is.

Beginning with (B), the comparison of both estimators indicates that $V_E\left(\hat{\Delta}\right) > V_A\left(\hat{\Delta}\right)$, except on two occasions, which means that the asymptotic formula itself always provides smaller values than those of the true variance. In particular, $V_E\left(\hat{\Delta}\right) - V_A\left(\hat{\Delta}\right)$ is found between $-0.0022$ and $+0.0310$ (with an average of 0.0063), which indicates that this underestimation may become very important.

Regarding what is highlighted in (A), the data always indicate that $\overline{\hat{V}\left(\hat{\Delta}\right)} > V_E\left(\hat{\Delta}\right)$ and $\overline{\hat{V}\left(\hat{\Delta}_U\right)} \geq V_E\left(\hat{\Delta}_U\right)$, so that (on average) the classic and new variances overestimate the real variance. It is also observed that $1.01 \leq \overline{\hat{V}\left(\hat{\Delta}\right)}/V_E\left(\hat{\Delta}\right) \leq 2.75$ and $1 \leq \overline{\hat{V}\left(\hat{\Delta}_U\right)}/V_E\left(\hat{\Delta}_U\right) \leq 2.71$, fractions which decrease (tending to 1) when $n$ or $K$ increase, but they grow when $\Delta$ increases. Note that although it was previously stated that $V_A\left(\hat{\Delta}\right)$ underestimated the real, value of the variance, substituting in it the correct optimal estimators the opposite occurs: the estimators overestimate the true variance. Note also that the formula used for $\hat{V}\left(\hat{\Delta}_U\right)$ works coherently in relation to the real formula of $\hat{V}\left(\hat{\Delta}\right)$.

*4.2. The case of parameter $\alpha_i$.*

The simulation procedure is similar to that of Section 4.1, leading to the results in Table 4 for the case of category 3 (the only one that is contemplated here, but the same thing happens with the other categories). Now the real value $\alpha_3$ is that of the setting itself, the $N$ estimators $\hat{\alpha}_3$ of each setting will give an average $\bar{\hat{\alpha}}_3$ and a variance $V_E\left(\hat{\alpha}_3\right)$ (denominator $N{-}1$) which is assumed to be the true one. We call $\overline{\hat{V}\left(\hat{\alpha}_3\right)}$ the average of the $N$ variances of $\bar{\hat{\alpha}}_3$. La exact asymptotic variance $V_A\left(\hat{\alpha}_3\right)$ is deduced from the values of the setting through the second expression (8). Similarly with $\hat{\alpha}_{3U}$ : $\bar{\hat{\alpha}}_{3U}$, $V_E\left(\hat{\alpha}_{3U}\right)$, and $\overline{\hat{V}\left(\hat{\alpha}_{3U}\right)}$.

It is observed that the two estimators always underestimate the true value $\alpha_3$, except for estimator $\hat{\alpha}_{3U}$ which they slightly overestimate on two occasions. The average underestimations are 0.017 for $\bar{\hat{\alpha}}_3$ and 0.008 for $\bar{\hat{\alpha}}_{3U}$. Furthermore $\bar{\hat{\alpha}}_{3U} \geq \bar{\hat{\alpha}}_3$, with a minimum of 0.001 and a maximum of 0.041. It is therefore clear that $\hat{\alpha}_{3U}$ is a better estimator of $\alpha_3$ than $\hat{\alpha}_3$. The average of the differences $\bar{\hat{\alpha}}_{3U} - \bar{\hat{\alpha}}_3$: (a) decreases with $n$, and is 0.014, 0.009, and 0.004 for $n = 30$, 50, and 100 respectively, which means that it is necessary to make the correction at least for $n{\leq}50$; (b) it decreases with $K$, and is 0.015 and 0.003 for $K{=}3$ and $K{=}5$ respectively, which means that it is necessary to make the correction at least for $K{\leq}3$; (c) it is not influenced excessively by $\alpha_3$. The consequence is that $\hat{\alpha}_{3U}$ is a better estimator than $\hat{\alpha}_3$ and it should be used at least for $n{\leq}50$ or $K{\leq}3$. For example, one of the extreme cases is that of the first setting ($n{=}30$, $K{=}3$, and $\alpha_3{=}0.20$) in which $\bar{\hat{\alpha}}_3{=}0.140$ and $\bar{\hat{\alpha}}_{3U}{=}0.182$. In a similar way as indicated in Section 4.1, in practice the demands may be rather less strict, since the values $\bar{\hat{\alpha}}_{3U} - \bar{\hat{\alpha}}_3$ are average differences and, therefore, the values of interest $\hat{\alpha}_{3U} - \hat{\alpha}_3$ may be much higher.

Regarding the variances, some similar happens to the scenario in the previous section: $V_E\left(\hat{\alpha}_3\right) > V_A\left(\hat{\alpha}_3\right)$, except on four occasions. This indicates that the asymptotic formula itself always provides values which are smaller than those of the true variance. Additionally, $V_E\left(\hat{\alpha}_3\right) - V_A\left(\hat{\alpha}_3\right)$ is found between $-0.0052$ and $+0.0321$ (with an average of 0.0034), which indicates that this underestimation may be very important. Furthermore, the data also always indicate that $\overline{\hat{V}\left(\hat{\alpha}_3\right)} \geq V_E\left(\hat{\alpha}_3\right)$ and $\overline{\hat{V}\left(\hat{\alpha}_{3U}\right)} \geq V_E\left(\hat{\alpha}_{3U}\right)$, so that (on average) the classic and new variances overestimate the real variance. It is also observed that $1 \leq \overline{\hat{V}\left(\hat{\alpha}_3\right)}/V_E\left(\hat{\alpha}_3\right) \leq 3.73$ and $1 \leq \overline{\hat{V}\left(\hat{\alpha}_{3U}\right)}/V_E\left(\hat{\alpha}_{3U}\right) \leq 2.63$, fractions that decrease (tending to 1) when $n$ or $K$ increase. Therefore, and as in the case of the previous section: (1) although it was previously stated that $V_A\left(\hat{\alpha}_3\right)$ underestimated the real value of the variance, substituting it with optimal estimators the opposite happens: the estimators overestimate the true variance; (2) the formula used for $\hat{V}\left(\hat{\alpha}_{iU}\right)$ works coherently in relation to the performance of the real formula real of $\hat{V}\left(\hat{\alpha}_i\right)$.

*4.3. The case of parameter $S_i$.*

The simulation procedure is similar to that in Section 4.1, leading to the results in Table 5 for the case of category 3. Now the real value $\mathcal{S}_3$ is deduced from the values of the setting in its row -see the last expression of (5)- and will be different in general for each row. The $N$=10,000 estimators $\hat{\mathcal{S}}_3$ of each setting will give a measure $\overline{\hat{\mathcal{S}}_3}$ and a variance $V_E\left(\hat{\mathcal{S}}_3\right)$ (denominator $N$–1) which is assumed to be the true one. We call $\overline{\hat{V}\left(\hat{\mathcal{S}}_3\right)}$ the average of the $N$ variances of $\hat{\mathcal{S}}_3$. The exact asymptotic variance $V_A\left(\hat{\mathcal{S}}_3\right)$ is deduced from the values of the setting through of the third expression (8). A similar thing happens with $\hat{\mathcal{S}}_{3U}$: $\overline{\hat{\mathcal{S}}_{3U}}$, $V_E\left(\hat{\mathcal{S}}_{3U}\right)$, and $\overline{\hat{V}\left(\hat{\mathcal{S}}_{3U}\right)}$.

It is observed that the two estimators always underestimate the true value $S_3$. The average underestimations are 0.05 for $\hat{\mathcal{S}}_3$ and 0.03 for $\hat{\mathcal{S}}_{3U}$. Moreover $\overline{\hat{\mathcal{S}}_{3U}} \geq \overline{\hat{\mathcal{S}}_3}$, with a minimum of 0.003 and a maximum of 0.094. It is therefore clear that $\hat{\mathcal{S}}_{3U}$ is clearly a better estimator of $\mathcal{S}_3$ than $\hat{\mathcal{S}}_3$. The average of the differences $\overline{\hat{\mathcal{S}}_{3U}} - \overline{\hat{\mathcal{S}}_3}$: (a) decreases with $n$, and is 0.04, 0.02, and 0.01 for $n = 30$, 50, and 100 respectively, which means that it is necessary to make the correction at least for $n \leq 50$; (b) it decreases with $K$, and is 0.04 and 0.01 for $K=3$ and $K=5$ respectively, which means that it is necessary to make the correction at least for $K \leq 3$; (c) it is not influenced excessively by $\mathcal{S}_3$. The consequence is that $\hat{\mathcal{S}}_{3U}$ is a better estimator than $\hat{\mathcal{S}}_3$ and it is necessary to use it at least for $n \leq 50$ or $K \leq 3$. For example, one of the extreme cases is that of the third setting ($n=30$, $K=3$, and $S_3=0.32$): $\overline{\hat{\mathcal{S}}_3}=0.19$ and $\overline{\hat{\mathcal{S}}_{3U}}=0.28$. In a similar way to what is indicated in Section 4.1, in practice the demands may be rather less strict, since the values $\overline{\hat{\mathcal{S}}_{3U}} - \overline{\hat{\mathcal{S}}_3}$ are average differences and, therefore, the values of interest $\hat{\mathcal{S}}_{3U} - \hat{\mathcal{S}}_3$ may be much higher.

It can be observed that it now occurs that $V_E\left(\hat{\mathcal{S}}_3\right) > V_A\left(\hat{\mathcal{S}}_3\right)$, except on two occasions, which indicates once again that the asymptotic formula itself always give smaller values than those of the true variance. In particular, $-0.0266 \leq V_E\left(\hat{\mathcal{S}}_3\right) - V_A\left(\hat{\mathcal{S}}_3\right) \leq +0.1608$, with an average of 0.0266; this indicates that this underestimation may also be very important. Furthermore, the data indicate that $\overline{V\left(\hat{\mathcal{S}}_3\right)} \geq V_E\left(\hat{\mathcal{S}}_3\right)$ and $\overline{V\left(\hat{\mathcal{S}}_{3U}\right)} \geq V_E\left(\hat{\mathcal{S}}_{3U}\right)$ on 75% of occasions, so that (on average) the classic and new variances usually overestimate the real variance. It is also observed that $0.87 \leq \overline{V\left(\hat{\mathcal{S}}_3\right)} / V_E\left(\hat{\mathcal{S}}_3\right) \leq 4.98$ and $0.96 \leq \overline{V\left(\hat{\mathcal{S}}_{3U}\right)} / V_E\left(\hat{\mathcal{S}}_{3U}\right) \leq 4.23$, fractions that decrease (tending to 1) when $n$ or $K$ increase. Note that it was previously stated that $V_A\left(\hat{\mathcal{S}}_3\right)$ almost always underestimated the real value of the variance $V_E\left(\hat{\mathcal{S}}_3\right)$, substituting in it the optimal estimators the opposite happens: the estimators usually overestimate the true variance. Additionally, the formula used for $\hat{V}\left(\hat{\mathcal{S}}_{3U}\right)$ works coherently in relation with the performance of the real formula of $\hat{V}\left(\hat{\mathcal{S}}_3\right)$.

## 5. The special case in which there are only two categories (*K*=2).

The *delta* model presents the problem of when there are only two categories, then there are more unknown parameters ($\alpha_1$, $\alpha_2$, $\pi_{11}$ and $\pi_{12}$) than cells free to take values ($p_{11}$, $p_{12}$, and $p_{21}$ for example, since $\Sigma_i\Sigma_j p_{ij}$=1). To solve the problem, Martín Andrés and Femia Marzo (2004) proposed the following solution: (1) create a third virtual category of observed frequencies $x_{i3}=x_{3j}=0$ ($\forall i, j$); (2) increase all of the data in the new 3×3 table by 0.5; (3) estimate the parameters based on the new table with *K*=3; and (4) redefine the measures of agreement taking into account the new situation. This procedure has been found to provide coherent results: Martín Andrés and Femia Marzo (2004, 2005, and 2008), Ato *et al.* (2011), Shankar and Bangdiwala (2014), and Giammarino *et al.* (2021).

Let $\overline{p}_{ij}$ be the new proportions observed, and let $\alpha_i$, $\pi_{ir}$ and $\Delta$ the parameters of the *delta* model, all of them referring to the new 3×3 table. With the notation of Martín Andrés and Álvarez Hernández (2022), the measures of agreement for the original 2×2 table are defined as $\alpha_i^*=\alpha_i/(p_{1\bullet}+p_{2\bullet})$, $\Delta^* = \alpha_1^* + \alpha_2^*$, and $\mathcal{S}_i^*=2\alpha_i/(p_{i\bullet}+p_{\bullet i})=\mathcal{S}_i$, for *i*=1 and 2. Their estimators and the estimated variances are:

$$\hat{\alpha}_i^* = \frac{\hat{\alpha}_i}{1-\overline{p}_{3\bullet}}, \quad \hat{V}\left(\hat{\alpha}_i^*\right) = \frac{1}{n\left(1-\overline{p}_{3\bullet}\right)^2}\left[\left(1-\hat{\Delta}\right)\hat{X}_i\left\{\frac{\hat{X}_i}{\hat{X}-1}-1\right\}+\left(1-\overline{p}_{3\bullet}\right)\hat{\alpha}_i^*\left(1-\hat{\alpha}_i^*\right)\right], \quad \textbf{(14)}$$

$$\hat{\Delta}^* = \hat{\alpha}_1^* + \hat{\alpha}_2^*, \quad \hat{V}\left(\hat{\Delta}^*\right) = \frac{1}{n\left(1-\overline{p}_{3\bullet}\right)^2}\left[\left(1-\hat{\Delta}\right)\left(1-\hat{X}_3\right)\frac{\hat{X}-\hat{X}_3}{\hat{X}-1}+\left(1-\overline{p}_{3\bullet}\right)\hat{\Delta}^*\left(1-\hat{\Delta}^*\right)\right]. \quad \textbf{(15)}$$

In the case of $\mathcal{S}_i$, whose definition is the same as that of expression (2), what is valid is what is indicated in expressions (5) for $\hat{\mathcal{S}}_i$ and (10) for $\hat{V}\left(\hat{\mathcal{S}}_i\right)$.

If now we apply to expressions (14) and (15) the corrections highlighted in Sections 2 and 3, we obtain the following expressions for the type $U$ estimators:

$$\hat{\alpha}_{iU}^* = \frac{\hat{\alpha}_{iU}}{1-\overline{p}_{3\bullet}}, \quad \hat{V}\left(\hat{\alpha}_{iU}^*\right) = \frac{1}{n\left(1-\overline{p}_{3\bullet}\right)^2}\left[\left(1-\hat{A}_U\right)\hat{X}_i\left\{\frac{\hat{X}_i}{\hat{X}-1}-1\right\} + \left(1-\overline{p}_{3\bullet}\right)\hat{\alpha}_{iU}^*\left(1-\hat{\alpha}_{iU}^*\right)\right], \quad \textbf{(16)}$$

$$\hat{A}_U^* = \hat{\alpha}_{1U}^* + \hat{\alpha}_{2U}^*, \quad \hat{V}\left(\hat{A}_U^*\right) = \frac{1}{n\left(1-\overline{p}_{3\bullet}\right)^2}\left[\left(1-\hat{A}_U\right)\left(1-\hat{X}_3\right)\frac{\hat{X}-\hat{X}_3}{\hat{X}-1} + \left(1-\overline{p}_{3\bullet}\right)\hat{A}_U^*\left(1-\hat{A}_U^*\right)\right]. \quad \textbf{(17)}$$

As in the previous paragraph, in the case of $\mathcal{S}_i$ expressions (7) for $\hat{\mathcal{S}}_{iU}$ and <span style="color:red">(12)</span> for $\hat{V}\left(\hat{\mathcal{S}}_{iU}\right)$ are valid.

## 6. Case in which there is a gold standard

Until now, it has been assumed that neither of two raters is a gold standard. Let us assume that the rater in rows is indeed a gold standard. In that case, two new parameters are also of interest. First, the raw probability that the rater will identify a subject as belonging to category $i$, when it actually belongs to that category, is $F_i = p_{ii}/p_{i\bullet}$; therefore that probability corrected by chance is the *conformity* $\mathcal{F}_i = \alpha_i/p_{i\bullet}$ by Martín Andrés and Femia Marzo (2005). Second, the raw probability that when the rater responses $i$, the subject is actually in the category $i$ (according to the gold standard), is $P_i = p_{ii}/p_{\bullet i}$; therefore that probability corrected by chance is the *predictivity* $\mathcal{P}_i = \alpha_i/p_{\bullet i}$ by the same authors. Note that when $K=2$, categories 1 and 2 are "sick" and "healthy" respectively, and the rater 2 is a diagnostic test, so $F_1$=Sensitivity, $F_2$=Specificity, $P_1$=Positive Predictive Value and $P_2$= Negative Predictive Value ("of the diagnostic test", in all cases). The estimators of these parameters, and their variances, are the ones by Martín Andrés and Femia Marzo (2005). Their values, in the format of Martín Andrés and Álvarez Hernández (2022), are:

$$\hat{\mathcal{F}}_i = \frac{\hat{\alpha}_i}{\overline{p}_{i\bullet}}, \ \hat{V}\left(\hat{\mathcal{F}}_i\right) = \frac{\hat{H}_i + \overline{p}_{i\bullet}\hat{\mathcal{F}}_i\left(1 - \hat{\mathcal{F}}_i\right)}{n\overline{p}_{i\bullet}^2}, \ \hat{\mathcal{P}}_i = \frac{\hat{\alpha}_i}{\overline{p}_{\bullet i}}, \ \text{and} \ \hat{V}\left(\hat{\mathcal{P}}_i\right) = \frac{\hat{H}_i + \overline{p}_{\bullet i}\hat{\mathcal{P}}_i\left(1 - \hat{\mathcal{P}}_i\right)}{n\overline{p}_{\bullet i}^2},$$

where $\hat{H}_i$ is given by the expression (11). Adapting those expressions to the current corrected estimators yields:

$$\hat{\mathcal{F}}_{iU} = \frac{\hat{\alpha}_{iU}}{\overline{p}_{i\bullet}}, \ \hat{V}\left(\hat{\mathcal{F}}_{iU}\right) = \frac{\hat{H}_{iU} + \overline{p}_{i\bullet}\hat{\mathcal{F}}_{iU}\left(1 - \hat{\mathcal{F}}_{iU}\right)}{n\overline{p}_{i\bullet}^2}, \ \hat{\mathcal{P}}_{iU} = \frac{\hat{\alpha}_{iU}}{\overline{p}_{\bullet i}}, \ \text{and} \ \hat{V}\left(\hat{\mathcal{P}}_{iU}\right) = \frac{\hat{H}_{iU} + \overline{p}_{\bullet i}\hat{\mathcal{P}}_{iU}\left(1 - \hat{\mathcal{P}}_{iU}\right)}{n\overline{p}_{\bullet i}^2},$$

where $\hat{H}_{iU}$ is given by the expression (13).

## 7. Examples

Table 6 (a) contains the data from the classic example by Fleiss *et al.* (2003) in which two raters diagnose $n$=100 subjects in $K$=3 categories (Psychotic, Neurotic and Organic). Part (b) specifies the classic estimations of parameters $\alpha_i$, $\mathcal{S}_i$, and $\Delta$ ($\hat{\alpha}_i$, $\hat{\mathcal{S}}_i$, and $\hat{\Delta}$, respectively) and the new estimations ($\hat{\alpha}_{iU}$, $\hat{\mathcal{S}}_{iU}$, and $\hat{\Delta}_U$, respectively). Some of the new estimators are higher than the classic estimators by approximately 0.03; for example, in the global agreement we obtain $\hat{\Delta}$= 0.687 *vs* $\hat{\Delta}_U$ = 0.715 (a 4% increase). With the classic Cohen *kappa* measure of agreement the differences are smaller: $\hat{\kappa}_C$ =0.676 *vs* $\hat{\kappa}_{CU}$ =0.679.

The differences are more notable when using the same $n$=100, but a $K$=2 value, as in Table 7(a); this table -which was already mentioned in the Introduction- also allow us to illustrate the advantage of the *delta* measure over the *kappa* measure. From the point of view of *kappa* it is obtained that $\hat{\kappa}_C$ =−0.111 *vs* $\hat{\kappa}_{CU}$ =−0.112. These values -which are very similar- are very surprising since, as the observed agreements are 80/100=80%, the two *kappa* coefficients provide a negative degree of agreement; this is due to the fact that the marginals are very unbalanced in the same direction. From the point of the view of *delta* it is obtained that -Table 7(b)- $\hat{\Delta}$= 0.582 *vs* $\hat{\Delta}_U$ = 0.714 (a 23% increase). Now the degree of agreement is notably larger, it is concordant with the high degree of agreement observed, and the value of $\hat{\Delta}_U$ is notably higher than the value of $\hat{\Delta}$. Something similar happens with the estimators of the coefficients $\alpha_i$ and $\mathcal{S}_i$. Furthermore, if we assume that the row rater is gold standard, the new interest coefficients that allow us evaluate the column rater are the *conformity* $\mathcal{F}_i$ and the *predictivity* $\mathcal{P}_i$. Their various estimators are those shown in the Table 7(c), and it can be seen that something similar to what is described above also occurs.

Finally, Table 8(a) contains the data from an example with *n*=30 subjects (*n*≤50) and *K*=4 categories (*K*>3). Part (b) specifies the classic estimations of parameters $\alpha_i$, $\mathcal{S}_i$ and $\Delta$ ($\hat{\alpha}_i$, $\hat{\mathcal{S}}_i$ and $\hat{\Delta}$, respectively) and the new estimations ($\hat{\alpha}_{iU}$, $\hat{\mathcal{S}}_{iU}$ and $\hat{\Delta}_U$, respectively). Some of the new estimators are higher than the classic estimators by 0.03 or more; for instance, in the global agreement we obtain $\hat{\Delta}$= 0.182 *vs* $\hat{\Delta}_U$ = 0.210 (a 15% increase) and in the consistency in category 2 we obtain $\hat{\mathcal{S}}_2$ =0.074 *vs* $\hat{\mathcal{S}}_{2U}$ =0.115 (a 55% increase). With the classic Cohen *kappa* measure of agreement the differences are smaller: $\hat{\kappa}_C$ =0.197 *vs* $\hat{\kappa}_{CU}$ =0.202.

## 8. Conclusions

There are different *kappa* type coefficients that measure the experimental degree of agreement between dos raters that independently classify *n* subjects in *K* categories. The best known of them is the $\kappa_C$ coefficient (Cohen 1960). Martín Andrés and Álvarez Hernández (2025) demonstrated that the traditional estimator $\hat{\kappa}_C$ of $\kappa_C$ underestimates the true value of the parameter, proposed a new estimator $\hat{\kappa}_{CU}$ ($\geq \hat{\kappa}_C$ si $\hat{\kappa}_C \geq 0$) that improves the performance of $\hat{\kappa}_C$, and justified that the correction is convenient at least when *n*≤30.

Given that $\kappa_C$ shows a poor performance when the marginal distributions are very unbalanced in the same direction (Brennan and Prediger 1981, Agresti *et al.* 1995, Guggenmoos-Holzmann and Vonk 1998, and Nelson and Pepe 2000) and, additionally, it does not always properly measure the degree of agreement in a specific category (Martín Andrés and Femia Marzo 2004, 2005), these same authors and Martín Andrés and Álvarez Hernández (2022) proposed a (*delta*) model que that corrected both problems. This article has demonstrated both theoretically and practically that the estimators $\hat{\alpha}_i, \hat{\mathcal{S}}_i$, and $\hat{\Delta}$ of the parameters of interest of the model $\alpha_i$, $\mathcal{S}_i$, and $\Delta$ respectively, usually underestimate the true value of said parameters, proposing new estimators $\hat{\alpha}_{iU}, \hat{\mathcal{S}}_{iU}$, and $\hat{\Delta}_U$ respectively, which improves the estimation. The corrected estimators should be used the smaller $n$ and $K$ are, and they are particularly advisable when $n \leq 50$ or $K \leq 3$.

The article also studies (through simulation) the behaviour of the variances of the previous estimators (classic and new). An initial conclusion is that the known formulas of the asymptotic variances of the estimators $\hat{\alpha}_i, \hat{\mathcal{S}}_i$, and $\hat{\Delta}$ underestimates the real variances of said parameters. Nevertheless, the estimations of said variances overestimates said real variances, especially in the case of the estimators $\hat{\alpha}_i$ and $\hat{\Delta}$; the same thing happens in the case of the new estimators $\hat{\alpha}_{iU}, \hat{\mathcal{S}}_{iU}$, and $\hat{\Delta}_U$.

Finally, the paper also considers two special cases: when $K=2$ and when the rater 1 is a gold standard. In the first case it is necessary to adapt formulas for the case $K>2$. In the second case, two new parameters are defined and estimated: the *conformity* $\mathcal{F}_i$ (capability of rater 2 to identify, not by chance, a subject belonging to the category $i$) and the *predictivity* $\mathcal{P}_i$ (reliability of a response $i$ from rater 2, not by chance). When both cases occur -$K=2$ and Rater 1 = Gold Standard-, the new parameters correspond to the correction by chance of classic parameters used to evaluate a diagnostic test: sensitivity and specificity on the one hand, and positive and negative predictive values on the other. In both cases, the new "$U$" type estimators perform better than the old uncorrected estimator of Martín Andrés and Femia Marzo (2005).

## Acknowledgments

This research was supported by the Ministry of Science and Innovation (Spain), Grant

PID2021-126095NB-I00 funded by MCIN/AEI/10.13039/501100011033 and by "ERDF A way of making Europe". The authors have declared no conflict of interest.

## APPENDIX A

**Variance and covariance of the estimators of $\pi_1$ and $\pi_2$ subject to the *delta* model**

Agresti (2013) particularized the multivariant delta model in the case of a multinomial distribution, pointing out that if the functions $f=f(\boldsymbol{p_{ij}})$ and $g=g(\boldsymbol{p_{ij}})$ are estimated by $\hat{f} = f\left(\boldsymbol{\bar{p}_{ij}}\right)$ and $\hat{g} = g\left(\boldsymbol{\bar{p}_{ij}}\right)$ with $\boldsymbol{p_{ij}}=\{p_{ij}\}$ and $\boldsymbol{\bar{p}_{ij}}=\left\{\bar{p}_{ij}\right\}$, then if $f_{ij}=\partial f/\partial p_{ij}$ and $g_{ij}=\partial g/\partial p_{ij}$:

$$nCov\left(\hat{f},\hat{g}\right) = \sum_{i=1}^{K}\sum_{j=1}^{K}f_{ij}g_{ij}p_{ij} - \left(\sum_{i=1}^{K}\sum_{j=1}^{K}f_{ij}\right)\left(\sum_{i=1}^{K}\sum_{j=1}^{K}g_{ij}\right). \tag{A1}$$

Moreover, Martín Andrés and Álvarez Hernández (2022) demonstrated that to obtain the estimators of the parameters of the *delta* model it is necessary to solve the $(K+1)$ equations that follow in the $(K+1)$ unknown $(\lambda_s, B)$, with $s=1, 2, \ldots, K$:

$$B = \frac{\left(\lambda_s + \bar{d}_{s1}\right)\left(\lambda_s + \bar{d}_{s2}\right)}{\lambda_s}, \quad \sum_{s=1}^{K}\lambda_s - B + \sum_{s=1}^{K}\bar{d}_{s1} = 0,$$

where $\bar{d}_{s1} = \bar{p}_{s\bullet} - \bar{p}_{ss}$ and $\bar{d}_{s2} = \bar{p}_{\bullet s} - \bar{p}_{ss}$ are, respectively, the observed disagreements of the raters 1 and 2 when classifying an individual like in category $s$. Note that $\sum_{s=1}^{K}\bar{d}_{s1}=\sum_{s=1}^{K}\bar{d}_{s2}$ is the total proportion of observed disagreements. In the previous expression it must be assumed that $\lambda_s=0$ when $\bar{d}_{ir}=0$ in $r=1$ or $r=2$. Once we have determined the values $\lambda_s$ and $B=1-\Delta$, the maximum likelihood estimators of the parameters of the model are:

$$\hat{\Delta} = 1 - B, \quad \hat{\alpha}_s = \bar{p}_{ss} - \lambda_s, \quad \hat{\pi}_{s1} = \frac{\lambda_s + \bar{d}_{s1}}{B}, \quad \hat{\pi}_{s2} = \frac{\lambda_s + \bar{d}_{s2}}{B} \quad \text{and} \quad \lambda_s = B\hat{\pi}_{s1}\hat{\pi}_{s2}. \tag{A2}$$

Note that if $\bar{d}_{s1}=0$, then $\hat{\alpha}_s = \bar{p}_{ss}$, $\hat{\pi}_{s1} = 0$, and $\hat{\pi}_{s2} = \bar{d}_{s2}/B$ (similarly when $\bar{d}_{s2}=0$). In that same article the following derivatives were also obtained,

$$\lambda_{B(s)} = \frac{\partial \lambda_s}{\partial B} = X_s, \quad \lambda_{ij(s)} = \frac{\partial \lambda_s}{\partial p_{ij}} = -\left(1-\delta_{ij}\right)X_s\left\{\frac{\delta_{si}}{\pi_{s1}} + \frac{\delta_{sj}}{\pi_{s2}}\right\},$$

$$\Delta_{ij} = \frac{\partial \Delta}{\partial p_{ij}} = \frac{1-\delta_{ij}}{X-1}\left\{1-\frac{X_i}{\pi_{i1}}-\frac{X_j}{\pi_{j2}}\right\} = -\frac{\partial B}{\partial p_{ij}}, \text{ and} \tag{A3}$$

$$\alpha_{ij(s)} = \frac{\partial \alpha_s}{\partial p_{ij}} = \delta_{si}\delta_{sj} + \left(1-\delta_{ij}\right)X_s\left[\left\{\frac{\delta_{si}}{\pi_{i1}}+\frac{\delta_{sj}}{\pi_{j2}}\right\}+\frac{1}{X-1}\left\{1-\frac{X_i}{\pi_{i1}}-\frac{X_j}{\pi_{j2}}\right\}\right],$$

with $X$ and $X_i$ as in Section 2. From this it is deduced that

$$\sum_s \lambda_{ij(s)} = -\left(1-\delta_{ij}\right)\left\{1-\frac{X_i}{\pi_{i1}}-\frac{X_j}{\pi_{j2}}\right\}, \text{ and}$$

$$\overline{\lambda}_{ij(s)} = \frac{d\lambda_s}{dp_{ij}} = -\left(1-\delta_{ij}\right)X_s\left[\left\{\frac{\delta_{si}}{\pi_{s1}}+\frac{\delta_{sj}}{\pi_{s2}}\right\}+\frac{1}{X-1}\left\{1-\frac{X_i}{\pi_{i1}}-\frac{X_j}{\pi_{j2}}\right\}\right], \tag{A4}$$

since $d\lambda_s/dp_{ij} = (\partial\lambda_s/\partial p_{ij})+(\partial\lambda_s/\partial\Delta)(\partial\Delta/\partial p_{ij}) = \lambda_{ij(s)}-X_s\Delta_{ij}$.

If $\hat{f}=\hat{\pi}_{s1}$ and $\hat{g}=\hat{\pi}_{s2}$ are the functions referred to at the beginning of this Appendix then, through expressions (A2), $f_{ij}=\partial\pi_{s1}/\partial p_{ij}=[B\{\overline{\lambda}_{ij(s)}+(\partial d_{s1}/\partial p_{ij})\}+(\lambda_s+d_{s1})\Delta_{ij}]/B^2$, with $d_{s1}=\Sigma_jp_{sj}-p_{ss}$ and $\partial d_{s1}/\partial p_{ij}=\delta_{si}(1-\delta_{ij})$. Substituting the second expression of (A3) and the expression (A4), taking into account that $\lambda_s+d_{s1}=B\pi_{s1}$ -third formula of expression (A2)-, and operating it is obtained that (in a similar way for $g_{ij}$):

$$f_{ij} = -\frac{1-\delta_{ij}}{(1-\Delta)T_s}\left[\delta_{si}\left(1-\pi_{s1}\right)+\delta_{sj}\pi_{s1}+\frac{\pi_{s1}\left(1-\pi_{s1}\right)}{X-1}\left\{1-\frac{X_i}{\pi_{i1}}-\frac{X_j}{\pi_{j2}}\right\}\right] \text{ and}$$

$$g_{ij} = -\frac{1-\delta_{ij}}{(1-\Delta)T_s}\left[\delta_{si}\pi_{s2}+\delta_{sj}\left(1-\pi_{s2}\right)+\frac{\pi_{s2}\left(1-\pi_{s2}\right)}{X-1}\left\{1-\frac{X_i}{\pi_{i1}}-\frac{X_j}{\pi_{j2}}\right\}\right].$$

where $T_s = \pi_{s1}+\pi_{s2}-1$. Performing operations it is deduced that $\Sigma_i\Sigma_jf_{ij}p_{ij} = \Sigma_i\Sigma_jg_{ij}p_{ij} = 0$ and that

$$\sum_i\sum_j f_{ij}g_{ij}p_{ij} = \frac{1}{1-\Delta}\left\{\pi_{s1}\pi_{s2}-\frac{X_s\left(X-X_s\right)}{X-1}\right\}, \tag{A5}$$

With this, applying expression (A1) it is deduced that $Cov\left(\hat{f},\hat{g}\right) = \left(\sum_{i=1}^K\sum_{j=1}^K f_{ij}g_{ij}p_{ij}\right)\!\big/n$, where $Cov\left(\hat{f},\hat{g}\right) = Cov\left(\hat{\pi}_{s1},\hat{\pi}_{s2}\right) = E\left(\hat{\pi}_{s1}\hat{\pi}_{s2}\right)-E\left(\hat{\pi}_{s1}\right)E\left(\hat{\pi}_{s2}\right)\approx E\left(\hat{\pi}_{s1}\hat{\pi}_{s2}\right)-\pi_{s1}\pi_{s2}$ since $\hat{\pi}_{sr}$ are maximum likelihood estimators. Therefore, for expression (A5), expression (3) is deduced from Section 2. Additionally:

$$V\left(\hat{\pi}_{sr}\right) = \frac{X_s - \pi_{sr}}{n\left(1-\Delta\right)}\left\{\frac{X_s - \pi_{sr}}{X-1} - 1\right\}$$

Section 2 proposed different estimators, all based on expression (4). An approximately unbiased estimator, which is an alternative to expression (4), is the following:

$$\left(\widehat{\pi_{i1}\pi_{i2}}\right)_{AC} = \frac{\hat{\pi}_{i1}\hat{\pi}_{i2} - \hat{A}_i}{\hat{C}} \quad \text{where} \quad \hat{A}_i = \frac{\hat{X}_i\left(\hat{X} - \hat{X}_i\right)}{n\left(1-\hat{\Delta}\right)\left(\hat{X}-1\right)} \quad \text{and} \quad \hat{C} = 1 + \frac{1}{n\left(1-\hat{\Delta}\right)},$$

which leads to the following alternative estimators to those of expressions (6) and (7):

$$\hat{I}_{\pi AC} = \frac{\hat{I}_\pi + A}{C}, \text{ where } \hat{A} = \sum_i \hat{A}_i = \frac{\hat{X}^2 - \sum_i \hat{X}_i^2}{n\left(1-\hat{\Delta}\right)\left(\hat{X}-1\right)},$$

$$\hat{\Delta}_{AC} = \frac{\hat{I}_0 - \hat{I}_{\pi AC}}{1 - \hat{I}_{\pi AC}}, \quad \hat{\alpha}_{iAC} = \overline{p}_{ii} - \left(1 - \hat{\Delta}_{AC}\right)\frac{\hat{\pi}_{i1}\hat{\pi}_{i2} + \hat{A}_i}{\hat{C}}, \text{ and } \hat{\mathcal{S}}_{iAC} = \frac{2\hat{\alpha}_{iAC}}{\overline{p}_{i\bullet} + \overline{p}_{\bullet i}}.$$

We have found that these estimators are only slightly worse than those reviewed in Section 2. Their approximate variances are the same as in expressions (12), changing $\hat{\Delta}_U$, $\hat{\alpha}_{iU}$, and $\hat{\mathcal{S}}_{iU}$ with $\hat{\Delta}_{AC}$, $\hat{\alpha}_{iAC}$, and $\hat{\mathcal{S}}_{iAC}$ respectively.

# Table 1

## Distribution of the responses when two observers classify $n$ subjects in $K$ categories

**(a) Absolute frequencies observed $x_{ij}$**

| Rater 2 → Rater 1 ↓ | 1 | . . . | $j$ | . . . | $K$ | Total |
|---|---|---|---|---|---|---|
| **1** | $x_{11}$ | . . . | $x_{1j}$ | . . . | $x_{1K}$ | $x_{1\bullet}$ |
| **.** | . | | . | | . | . |
| **.** | . | | . | | . | . |
| **.** | . | | . | | . | . |
| ***i*** | $x_{i1}$ | . . . | $x_{ij}$ | . . . | $x_{iK}$ | $x_{i\bullet}$ |
| **.** | . | | . | | . | . |
| **.** | . | | . | | . | . |
| ***K*** | $x_{K1}$ | . . . | $x_{Kj}$ | . . . | $x_{KK}$ | $x_{K\bullet}$ |
| **Total** | $x_{\bullet1}$ | . . . | $x_{\bullet j}$ | . . . | $x_{\bullet K}$ | $n$ |

**(b) Probability of a subject being classified in the indicated cell**

| Rater 2 → Rater 1 ↓ | 1 | . . . | $J$ | . . . | $K$ | Total |
|---|---|---|---|---|---|---|
| **1** | $p_{11}$ | . . . | $p_{1j}$ | . . . | $p_{1K}$ | $p_{1\bullet}$ |
| **.** | . | | . | | . | . |
| **.** | . | | . | | . | . |
| **.** | . | | . | | . | . |
| ***i*** | $p_{i1}$ | . . . | $p_{ij}$ | . . . | $p_{iK}$ | $p_{i\bullet}$ |
| **.** | . | | . | | . | . |
| **.** | . | | . | | . | . |
| ***K*** | $p_{K1}$ | . . . | $p_{Kj}$ | . . . | $p_{KK}$ | $p_{K\bullet}$ |
| **Total** | $p_{\bullet1}$ | . . . | $p_{\bullet j}$ | . . . | $p_{\bullet K}$ | 1 |

**Table 2: Each line of the table indicates the populational values to be used (values of _K_, _n_, _αᵢ_ and _πᵢᵣ_). The first column (id) refers to the number of the setting.**

| id | K | n | $\alpha_1$ | $\alpha_2$ | $\alpha_3$ | $\alpha_4$ | $\alpha_5$ | $\pi_{11}$ | $\pi_{21}$ | $\pi_{31}$ | $\pi_{41}$ | $\pi_{51}$ | $\pi_{12}$ | $\pi_{22}$ | $\pi_{32}$ | $\pi_{42}$ | $\pi_{52}$ |
|---|---|---|---|---|---|---|---|---|---|---|---|---|---|---|---|---|---|
| 1 | 3 | 30 | 0.05 | 0.15 | 0.2 | --- | --- | 0.2 | 0.3 | 0.5 | --- | --- | 0.2 | 0.3 | 0.5 | --- | --- |
| 2 | | | | | | --- | --- | | | | --- | --- | 0.5 | 0.3 | 0.2 | --- | --- |
| 3 | | | 0.13 | 0.13 | 0.14 | --- | --- | | | | --- | --- | 0.2 | 0.3 | 0.5 | --- | --- |
| 4 | | | | | | --- | --- | | | | --- | --- | 0.5 | 0.3 | 0.2 | --- | --- |
| 5 | | 50 | 0.05 | 0.15 | 0.2 | --- | --- | | | | --- | --- | 0.2 | 0.3 | 0.5 | --- | --- |
| 6 | | | | | | --- | --- | | | | --- | --- | 0.5 | 0.3 | 0.2 | --- | --- |
| 7 | | | 0.13 | 0.13 | 0.14 | --- | --- | | | | --- | --- | 0.2 | 0.3 | 0.5 | --- | --- |
| 8 | | | | | | --- | --- | | | | --- | --- | 0.5 | 0.3 | 0.2 | --- | --- |
| 9 | | 100 | 0.05 | 0.15 | 0.2 | --- | --- | | | | --- | --- | 0.2 | 0.3 | 0.5 | --- | --- |
| 10 | | | | | | --- | --- | | | | --- | --- | 0.5 | 0.3 | 0.2 | --- | --- |
| 11 | | | 0.13 | 0.13 | 0.14 | --- | --- | | | | --- | --- | 0.2 | 0.3 | 0.5 | --- | --- |
| 12 | | | | | | --- | --- | | | | --- | --- | 0.5 | 0.3 | 0.2 | --- | --- |
| 13 | | 30 | 0.15 | 0.25 | 0.4 | --- | --- | | | | --- | --- | 0.2 | 0.3 | 0.5 | --- | --- |
| 14 | | | | | | --- | --- | | | | --- | --- | 0.5 | 0.3 | 0.2 | --- | --- |
| 15 | | | 0.26 | 0.26 | 0.28 | --- | --- | | | | --- | --- | 0.2 | 0.3 | 0.5 | --- | --- |
| 16 | | | | | | --- | --- | | | | --- | --- | 0.5 | 0.3 | 0.2 | --- | --- |
| 17 | | 50 | 0.15 | 0.25 | 0.4 | --- | --- | | | | --- | --- | 0.2 | 0.3 | 0.5 | --- | --- |
| 18 | | | | | | --- | --- | | | | --- | --- | 0.5 | 0.3 | 0.2 | --- | --- |
| 19 | | | 0.26 | 0.26 | 0.28 | --- | --- | | | | --- | --- | 0.2 | 0.3 | 0.5 | --- | --- |
| 20 | | | | | | --- | --- | | | | --- | --- | 0.5 | 0.3 | 0.2 | --- | --- |
| 21 | | 100 | 0.15 | 0.25 | 0.4 | --- | --- | | | | --- | --- | 0.2 | 0.3 | 0.5 | --- | --- |
| 22 | | | | | | --- | --- | | | | --- | --- | 0.5 | 0.3 | 0.2 | --- | --- |
| 23 | | | 0.26 | 0.26 | 0.28 | --- | --- | | | | --- | --- | 0.2 | 0.3 | 0.5 | --- | --- |
| 24 | | | | | | --- | --- | | | | --- | --- | 0.5 | 0.3 | 0.2 | --- | --- |
| 25 | 5 | 30 | 0.05 | 0.05 | 0.05 | 0.1 | 0.15 | 0.1 | 0.15 | 0.2 | 0.25 | 0.3 | 0.1 | 0.15 | 0.2 | 0.25 | 0.3 |
| 26 | | | | | | | | | | | | | 0.3 | 0.25 | 0.2 | 0.15 | 0.1 |
| 27 | | | 0.08 | 0.08 | 0.08 | 0.08 | 0.08 | | | | | | 0.1 | 0.15 | 0.2 | 0.25 | 0.3 |
| 28 | | | | | | | | | | | | | 0.3 | 0.25 | 0.2 | 0.15 | 0.1 |
| 28 | | 50 | 0.05 | 0.05 | 0.05 | 0.1 | 0.15 | | | | | | 0.1 | 0.15 | 0.2 | 0.25 | 0.3 |
| 30 | | | | | | | | | | | | | 0.3 | 0.25 | 0.2 | 0.15 | 0.1 |
| 31 | | | 0.08 | 0.08 | 0.08 | 0.08 | 0.08 | | | | | | 0.1 | 0.15 | 0.2 | 0.25 | 0.3 |
| 32 | | | | | | | | | | | | | 0.3 | 0.25 | 0.2 | 0.15 | 0.1 |
| 33 | | 100 | 0.05 | 0.05 | 0.05 | 0.1 | 0.15 | | | | | | 0.1 | 0.15 | 0.2 | 0.25 | 0.3 |
| 34 | | | | | | | | | | | | | 0.3 | 0.25 | 0.2 | 0.15 | 0.1 |
| 35 | | | 0.08 | 0.08 | 0.08 | 0.08 | 0.08 | | | | | | 0.1 | 0.15 | 0.2 | 0.25 | 0.3 |
| 36 | | | | | | | | | | | | | 0.3 | 0.25 | 0.2 | 0.15 | 0.1 |
| 37 | | 30 | 0.1 | 0.15 | 0.15 | 0.2 | 0.2 | | | | | | 0.1 | 0.15 | 0.2 | 0.25 | 0.3 |
| 38 | | | | | | | | | | | | | 0.3 | 0.25 | 0.2 | 0.15 | 0.1 |
| 38 | | | 0.16 | 0.16 | 0.16 | 0.16 | 0.16 | | | | | | 0.1 | 0.15 | 0.2 | 0.25 | 0.3 |
| 40 | | | | | | | | | | | | | 0.3 | 0.25 | 0.2 | 0.15 | 0.1 |
| 41 | | 50 | 0.1 | 0.15 | 0.15 | 0.2 | 0.2 | | | | | | 0.1 | 0.15 | 0.2 | 0.25 | 0.3 |
| 42 | | | | | | | | | | | | | 0.3 | 0.25 | 0.2 | 0.15 | 0.1 |
| 43 | | | 0.16 | 0.16 | 0.16 | 0.16 | 0.16 | | | | | | 0.1 | 0.15 | 0.2 | 0.25 | 0.3 |
| 44 | | | | | | | | | | | | | 0.3 | 0.25 | 0.2 | 0.15 | 0.1 |
| 45 | | 100 | 0.1 | 0.15 | 0.15 | 0.2 | 0.2 | | | | | | 0.1 | 0.15 | 0.2 | 0.25 | 0.3 |
| 46 | | | | | | | | | | | | | 0.3 | 0.25 | 0.2 | 0.15 | 0.1 |
| 47 | | | 0.16 | 0.16 | 0.16 | 0.16 | 0.16 | | | | | | 0.1 | 0.15 | 0.2 | 0.25 | 0.3 |
| 48 | | | | | | | | | | | | | 0.3 | 0.25 | 0.2 | 0.15 | 0.1 |

**Table 3: Assessment of the estimators $\hat{\Delta}$ and $\hat{\Delta}_U$, and of their variances, based on 10,000 simulation of each one of the settings in Table 2.**

| id | $K$ | $n$ | $\Delta$ | $\overline{\hat{\Delta}}$ | $\overline{\hat{\Delta}_U}$ | $V_A(\hat{\Delta})$ | $V_E(\hat{\Delta})$ | $\overline{\hat{V}(\hat{\Delta})}$ | $V_E(\hat{\Delta}_U)$ | $\overline{\hat{V}(\hat{\Delta}_U)}$ |
|----|----|-----|-----|--------|--------|--------|--------|--------|--------|--------|
| 1 | 3 | 30 | 0.4 | 0.3127 | 0.3824 | 0.0280 | 0.0509 | 0.1214 | 0.0413 | 0.1091 |
| 2 | | | | 0.3416 | 0.3829 | 0.0174 | 0.0321 | 0.0641 | 0.0274 | 0.0596 |
| 3 | | | | 0.3171 | 0.3867 | 0.0280 | 0.0485 | 0.1179 | 0.0395 | 0.1059 |
| 4 | | | | 0.3373 | 0.3786 | 0.0174 | 0.0312 | 0.0653 | 0.0266 | 0.0607 |
| 5 | | 50 | | 0.3317 | 0.3752 | 0.0168 | 0.0404 | 0.0896 | 0.0330 | 0.0822 |
| 6 | | | | 0.3692 | 0.3905 | 0.0105 | 0.0161 | 0.0248 | 0.0147 | 0.0238 |
| 7 | | | | 0.3238 | 0.3678 | 0.0168 | 0.0448 | 0.0976 | 0.0367 | 0.0896 |
| 8 | | | | 0.3675 | 0.3892 | 0.0105 | 0.0164 | 0.0262 | 0.0149 | 0.0252 |
| 9 | | 100 | | 0.3664 | 0.3854 | 0.0084 | 0.0206 | 0.0302 | 0.0182 | 0.0288 |
| 10 | | | | 0.3873 | 0.3965 | 0.0052 | 0.0060 | 0.0066 | 0.0059 | 0.0066 |
| 11 | | | | 0.3647 | 0.3838 | 0.0084 | 0.0204 | 0.0299 | 0.0181 | 0.0286 |
| 12 | | | | 0.3892 | 0.3984 | 0.0052 | 0.0057 | 0.0065 | 0.0056 | 0.0065 |
| 13 | | 30 | 0.8 | 0.7088 | 0.7513 | 0.0120 | 0.0101 | 0.0224 | 0.0077 | 0.0196 |
| 14 | | | | 0.7083 | 0.7446 | 0.0085 | 0.0096 | 0.0223 | 0.0074 | 0.0199 |
| 15 | | | | 0.7095 | 0.7521 | 0.0120 | 0.0098 | 0.0223 | 0.0075 | 0.0195 |
| 16 | | | | 0.7079 | 0.7442 | 0.0085 | 0.0093 | 0.0222 | 0.0073 | 0.0198 |
| 17 | | 50 | | 0.7479 | 0.7815 | 0.0072 | 0.0091 | 0.0189 | 0.0071 | 0.0166 |
| 18 | | | | 0.7543 | 0.7788 | 0.0051 | 0.0077 | 0.0156 | 0.0063 | 0.0140 |
| 19 | | | | 0.7497 | 0.7833 | 0.0072 | 0.0093 | 0.0189 | 0.0073 | 0.0166 |
| 20 | | | | 0.7540 | 0.7786 | 0.0051 | 0.0076 | 0.0159 | 0.0062 | 0.0143 |
| 21 | | 100 | | 0.7708 | 0.7923 | 0.0036 | 0.0067 | 0.0145 | 0.0055 | 0.0130 |
| 22 | | | | 0.7816 | 0.7937 | 0.0025 | 0.0039 | 0.0071 | 0.0034 | 0.0066 |
| 23 | | | | 0.7698 | 0.7913 | 0.0036 | 0.0070 | 0.0150 | 0.0058 | 0.0134 |
| 24 | | | | 0.7809 | 0.7931 | 0.0025 | 0.0040 | 0.0073 | 0.0035 | 0.0068 |
| 25 | 5 | 30 | 0.4 | 0.3629 | 0.3802 | 0.0143 | 0.0186 | 0.0227 | 0.0173 | 0.0221 |
| 26 | | | | 0.3797 | 0.3908 | 0.0125 | 0.0140 | 0.0152 | 0.0137 | 0.0150 |
| 27 | | | | 0.3699 | 0.3873 | 0.0143 | 0.0182 | 0.0233 | 0.0171 | 0.0227 |
| 28 | | | | 0.3824 | 0.3934 | 0.0125 | 0.0140 | 0.0151 | 0.0137 | 0.0150 |
| 28 | | 50 | | 0.3873 | 0.3967 | 0.0086 | 0.0094 | 0.0103 | 0.0092 | 0.0102 |
| 30 | | | | 0.3934 | 0.3995 | 0.0075 | 0.0076 | 0.0081 | 0.0076 | 0.0080 |
| 31 | | | | 0.3875 | 0.3970 | 0.0086 | 0.0095 | 0.0103 | 0.0093 | 0.0102 |
| 32 | | | | 0.3919 | 0.3980 | 0.0075 | 0.0080 | 0.0081 | 0.0079 | 0.0081 |
| 33 | | 100 | | 0.3952 | 0.3997 | 0.0043 | 0.0045 | 0.0046 | 0.0045 | 0.0045 |
| 34 | | | | 0.3982 | 0.4010 | 0.0038 | 0.0038 | 0.0039 | 0.0038 | 0.0039 |
| 35 | | | | 0.3941 | 0.3985 | 0.0043 | 0.0045 | 0.0046 | 0.0045 | 0.0045 |
| 36 | | | | 0.3971 | 0.4000 | 0.0038 | 0.0038 | 0.0039 | 0.0038 | 0.0039 |
| 37 | | 30 | 0.8 | 0.6347 | 0.6912 | 0.0074 | 0.0384 | 0.0993 | 0.0173 | 0.0414 |
| 38 | | | | 0.6836 | 0.7175 | 0.0068 | 0.0266 | 0.0549 | 0.0133 | 0.0269 |
| 38 | | | | 0.6386 | 0.6968 | 0.0074 | 0.0376 | 0.1036 | 0.0165 | 0.0426 |
| 40 | | | | 0.6885 | 0.7215 | 0.0068 | 0.0253 | 0.0534 | 0.0125 | 0.0260 |
| 41 | | 50 | | 0.7650 | 0.7816 | 0.0045 | 0.0087 | 0.0163 | 0.0063 | 0.0107 |
| 42 | | | | 0.7780 | 0.7877 | 0.0041 | 0.0061 | 0.0088 | 0.0050 | 0.0072 |
| 43 | | | | 0.7630 | 0.7797 | 0.0045 | 0.0085 | 0.0163 | 0.0062 | 0.0105 |
| 44 | | | | 0.7784 | 0.7880 | 0.0041 | 0.0058 | 0.0088 | 0.0048 | 0.0071 |
| 45 | | 100 | | 0.7933 | 0.7983 | 0.0022 | 0.0025 | 0.0029 | 0.0024 | 0.0028 |
| 46 | | | | 0.7958 | 0.7991 | 0.0021 | 0.0021 | 0.0023 | 0.0021 | 0.0023 |
| 47 | | | | 0.7935 | 0.7986 | 0.0022 | 0.0024 | 0.0030 | 0.0024 | 0.0029 |
| 48 | | | | 0.7960 | 0.7993 | 0.0021 | 0.0021 | 0.0023 | 0.0021 | 0.0023 |

(1) $\Delta = \Sigma_i \alpha_i$ is the true value of the *delta* parameter; $\overline{\hat{\Delta}}$ and $\overline{\hat{\Delta}_U}$ are the sample averages of the 10,000 estimators $\hat{\Delta}$ and $\hat{\Delta}_U$ respectively.

(2) $V_A(\hat{\Delta})$ is the real asymptotic variance, obtained based on the true values of the parameters.

(3) $V_E\left(\hat{\Delta}\right)$ and $V_E\left(\hat{\Delta}_U\right)$ are the sample variances of the 10,000 estimators $\hat{\Delta}$ and $\hat{\Delta}_U$ respectively. It is assumed that these variances are approximately equal to the true variances of the estimators $\hat{\Delta}$ and $\hat{\Delta}_U$ respectively.

(4) $\overline{\hat{V}\left(\hat{\Delta}\right)}$ and $\overline{\hat{V}\left(\hat{\Delta}_U\right)}$ are the sample averages of the 10,000 estimators $\hat{V}\left(\hat{\Delta}\right)$ and $\hat{V}\left(\hat{\Delta}_U\right)$ respectively.

**Table 4: Assessment of the estimators of $\alpha_3$ ($\hat{\alpha}_3$ and $\hat{\alpha}_{3E}$), and of their variances, based on 10,000 simulations of each one of the settings in Table 2.**

| id | K | n | $\alpha_3$ | $\bar{\bar{\alpha}}_3$ | $\bar{\bar{\alpha}}_{3U}$ | $V_A(\hat{\alpha}_3)$ | $V_E(\hat{\alpha}_3)$ | $V_E(\hat{\alpha}_{3U})$ | $\overline{\hat{V}(\hat{\alpha}_3)}$ | $\overline{\hat{V}(\hat{\alpha}_{3U})}$ |
|---|---|---|---|---|---|---|---|---|---|---|
| 1 | 3 | 30 | 0.2 | 0.1402 | 0.1815 | 0.0312 | 0.0518 | 0.0406 | 0.1133 | 0.1018 |
| 2 | | | | 0.1785 | 0.1930 | 0.0108 | 0.0196 | 0.0164 | 0.0305 | 0.0286 |
| 3 | | | 0.14 | 0.0859 | 0.1270 | 0.0298 | 0.0476 | 0.0371 | 0.1085 | 0.0975 |
| 4 | | | | 0.1146 | 0.1296 | 0.0095 | 0.0186 | 0.0152 | 0.0312 | 0.0291 |
| 5 | | 50 | 0.2 | 0.1469 | 0.1742 | 0.0187 | 0.0450 | 0.0365 | 0.0875 | 0.0803 |
| 6 | | | | 0.1865 | 0.1942 | 0.0065 | 0.0105 | 0.0094 | 0.0139 | 0.0134 |
| 7 | | | 0.14 | 0.0783 | 0.1063 | 0.0179 | 0.0500 | 0.0404 | 0.0961 | 0.0883 |
| 8 | | | | 0.1285 | 0.1361 | 0.0057 | 0.0093 | 0.0083 | 0.0123 | 0.0119 |
| 9 | | 100 | 0.2 | 0.1709 | 0.1829 | 0.0094 | 0.0236 | 0.0209 | 0.0314 | 0.0300 |
| 10 | | | | 0.1953 | 0.1986 | 0.0032 | 0.0040 | 0.0038 | 0.0040 | 0.0040 |
| 11 | | | 0.14 | 0.1080 | 0.1201 | 0.0090 | 0.0230 | 0.0203 | 0.0308 | 0.0294 |
| 12 | | | | 0.1349 | 0.1382 | 0.0029 | 0.0035 | 0.0033 | 0.0036 | 0.0036 |
| 13 | | 30 | 0.4 | 0.3586 | 0.3786 | 0.0166 | 0.0114 | 0.0097 | 0.0201 | 0.0184 |
| 14 | | | | 0.3582 | 0.3707 | 0.0098 | 0.0095 | 0.0085 | 0.0142 | 0.0133 |
| 15 | | | 0.28 | 0.2516 | 0.2716 | 0.0153 | 0.0104 | 0.0088 | 0.0187 | 0.0171 |
| 16 | | | | 0.2480 | 0.2607 | 0.0086 | 0.0085 | 0.0075 | 0.0131 | 0.0123 |
| 17 | | 50 | 0.4 | 0.3742 | 0.3918 | 0.0100 | 0.0101 | 0.0085 | 0.0181 | 0.0163 |
| 18 | | | | 0.3808 | 0.3893 | 0.0059 | 0.0071 | 0.0063 | 0.0099 | 0.0093 |
| 19 | | | 0.28 | 0.2615 | 0.2791 | 0.0092 | 0.0094 | 0.0078 | 0.0172 | 0.0155 |
| 20 | | | | 0.2645 | 0.2732 | 0.0051 | 0.0064 | 0.0056 | 0.0093 | 0.0088 |
| 21 | | 100 | 0.4 | 0.3801 | 0.3930 | 0.0050 | 0.0080 | 0.0067 | 0.0148 | 0.0134 |
| 22 | | | | 0.3925 | 0.3968 | 0.0029 | 0.0038 | 0.0035 | 0.0049 | 0.0047 |
| 23 | | | 0.28 | 0.2597 | 0.2727 | 0.0046 | 0.0077 | 0.0064 | 0.0147 | 0.0134 |
| 24 | | | | 0.2727 | 0.2771 | 0.0026 | 0.0034 | 0.0031 | 0.0046 | 0.0044 |
| 25 | 5 | 30 | 0.05 | 0.0440 | 0.0471 | 0.0029 | 0.0034 | 0.0032 | 0.0038 | 0.0038 |
| 26 | | | | 0.0462 | 0.0487 | 0.0028 | 0.0033 | 0.0031 | 0.0036 | 0.0035 |
| 27 | | | 0.08 | 0.0737 | 0.0768 | 0.0037 | 0.0042 | 0.0040 | 0.0045 | 0.0045 |
| 28 | | | | 0.0761 | 0.0785 | 0.0037 | 0.0042 | 0.0040 | 0.0043 | 0.0044 |
| 28 | | 50 | 0.05 | 0.0477 | 0.0494 | 0.0017 | 0.0019 | 0.0018 | 0.0019 | 0.0019 |
| 30 | | | | 0.0483 | 0.0497 | 0.0017 | 0.0018 | 0.0018 | 0.0019 | 0.0019 |
| 31 | | | 0.08 | 0.0795 | 0.0811 | 0.0022 | 0.0023 | 0.0023 | 0.0024 | 0.0024 |
| 32 | | | | 0.0784 | 0.0798 | 0.0022 | 0.0024 | 0.0024 | 0.0024 | 0.0024 |
| 33 | | 100 | 0.05 | 0.0493 | 0.0500 | 0.0009 | 0.0009 | 0.0009 | 0.0009 | 0.0009 |
| 34 | | | | 0.0497 | 0.0503 | 0.0009 | 0.0009 | 0.0009 | 0.0009 | 0.0009 |
| 35 | | | 0.08 | 0.0788 | 0.0796 | 0.0011 | 0.0011 | 0.0011 | 0.0011 | 0.0011 |
| 36 | | | | 0.0794 | 0.0801 | 0.0011 | 0.0011 | 0.0011 | 0.0011 | 0.0011 |
| 37 | | 30 | 0.15 | 0.1199 | 0.1305 | 0.0047 | 0.0071 | 0.0047 | 0.0262 | 0.0118 |
| 38 | | | | 0.1277 | 0.1347 | 0.0047 | 0.0066 | 0.0048 | 0.0175 | 0.0095 |
| 38 | | | 0.16 | 0.1289 | 0.1397 | 0.0049 | 0.0074 | 0.0049 | 0.0276 | 0.0122 |
| 40 | | | | 0.1358 | 0.1428 | 0.0049 | 0.0066 | 0.0049 | 0.0180 | 0.0096 |
| 41 | | 50 | 0.15 | 0.1437 | 0.1467 | 0.0028 | 0.0032 | 0.0029 | 0.0050 | 0.0038 |
| 42 | | | | 0.1462 | 0.1483 | 0.0028 | 0.0032 | 0.0030 | 0.0037 | 0.0035 |
| 43 | | | 0.16 | 0.1527 | 0.1556 | 0.0029 | 0.0032 | 0.0030 | 0.0046 | 0.0038 |
| 44 | | | | 0.1557 | 0.1578 | 0.0029 | 0.0032 | 0.0031 | 0.0040 | 0.0036 |
| 45 | | 100 | 0.15 | 0.1491 | 0.1499 | 0.0014 | 0.0015 | 0.0014 | 0.0015 | 0.0015 |
| 46 | | | | 0.1500 | 0.1507 | 0.0014 | 0.0014 | 0.0014 | 0.0015 | 0.0015 |
| 47 | | | 0.16 | 0.1592 | 0.1601 | 0.0015 | 0.0015 | 0.0015 | 0.0015 | 0.0015 |
| 48 | | | | 0.1588 | 0.1596 | 0.0015 | 0.0015 | 0.0015 | 0.0015 | 0.0015 |

(1) $\alpha_3$ is the true value of the agreement in category 3; $\bar{\bar{\alpha}}_3$ and $\bar{\bar{\alpha}}_{3U}$ are the sample averages of the 10,000 estimators $\hat{\alpha}_3$ and $\hat{\alpha}_{3U}$ respectively.

(2) $V_A(\hat{\alpha}_3)$ is the real asymptotic variance, obtained based on the true values of the parameters.

(3) $V_E(\hat{\alpha}_3)$ and $V_E(\hat{\alpha}_{3U})$ are the sample variances of the 10,000 estimators $\hat{\alpha}_3$ and $\hat{\alpha}_{3U}$ respectively. It is assumed that those sample variances are approximately equal to the true variances of the previous estimators.

(4) $\overline{\hat{V}(\hat{\alpha}_3)}$ and $\overline{\hat{V}(\hat{\alpha}_{3U})}$ are the sample averages of the 10,000 estimators $\hat{V}(\hat{\alpha}_3)$ and $\hat{V}(\hat{\alpha}_{3U})$ respectively.

**Table 5: Evaluation of three estimators of $\mathcal{S}_3$ ($\hat{\mathcal{S}}_3$ and $\hat{\mathcal{S}}_{3U}$), and of their variances, based on 10,000 simulations of each one of the settings in Table 2.**

| id | $K$ | $n$ | $\mathcal{S}_3$ | $\bar{\hat{\mathcal{S}}}_3$ | $\bar{\hat{\mathcal{S}}}_{3U}$ | $V_A\left(\hat{\mathcal{S}}_3\right)$ | $V_E\left(\hat{\mathcal{S}}_3\right)$ | $V_E\left(\hat{\mathcal{S}}_{3U}\right)$ | $\overline{\hat{V}\left(\hat{\mathcal{S}}_3\right)}$ | $\overline{\hat{V}\left(\hat{\mathcal{S}}_{3U}\right)}$ |
|----|----|----|----|----|----|----|----|----|----|----|
| 1 | 3 | 30 | 0.4000 | 0.2756 | 0.3590 | 0.1177 | 0.2093 | 0.1602 | 0.4471 | 0.4012 |
| 2 | | | 0.4878 | 0.4296 | 0.4645 | 0.0494 | 0.0999 | 0.0803 | 0.1512 | 0.1408 |
| 3 | | | 0.3182 | 0.1896 | 0.2833 | 0.1486 | 0.2483 | 0.1897 | 0.5400 | 0.4857 |
| 4 | | | 0.4000 | 0.3228 | 0.3647 | 0.0644 | 0.1335 | 0.1061 | 0.2092 | 0.1954 |
| 5 | | 50 | 0.4000 | 0.2930 | 0.3470 | 0.0706 | 0.1755 | 0.1410 | 0.3293 | 0.3021 |
| 6 | | | 0.4878 | 0.4521 | 0.4706 | 0.0297 | 0.0523 | 0.0456 | 0.0665 | 0.0637 |
| 7 | | | 0.3182 | 0.1764 | 0.2394 | 0.0891 | 0.2499 | 0.2008 | 0.4632 | 0.4256 |
| 8 | | | 0.4000 | 0.3632 | 0.3846 | 0.0387 | 0.0651 | 0.0569 | 0.0818 | 0.0787 |
| 9 | | 100 | 0.4000 | 0.3416 | 0.3653 | 0.0353 | 0.0899 | 0.0793 | 0.1167 | 0.1115 |
| 10 | | | 0.4878 | 0.4741 | 0.4819 | 0.0148 | 0.0189 | 0.0179 | 0.0190 | 0.0187 |
| 11 | | | 0.3182 | 0.2458 | 0.2730 | 0.0446 | 0.1122 | 0.0988 | 0.1461 | 0.1395 |
| 12 | | | 0.4000 | 0.3836 | 0.3929 | 0.0193 | 0.0243 | 0.0230 | 0.0244 | 0.0242 |
| 13 | | 30 | 0.8000 | 0.7378 | 0.7797 | 0.0430 | 0.0283 | 0.0196 | 0.0684 | 0.0587 |
| 14 | | | 0.8511 | 0.7760 | 0.8036 | 0.0145 | 0.0204 | 0.0143 | 0.0442 | 0.0386 |
| 15 | | | 0.7368 | 0.6608 | 0.7153 | 0.0735 | 0.0469 | 0.0325 | 0.1118 | 0.0965 |
| 16 | | | 0.8000 | 0.6995 | 0.7366 | 0.0258 | 0.0356 | 0.0250 | 0.0769 | 0.0675 |
| 17 | | 50 | 0.8000 | 0.7556 | 0.7916 | 0.0258 | 0.0286 | 0.0208 | 0.0627 | 0.0544 |
| 18 | | | 0.8511 | 0.8140 | 0.8324 | 0.0087 | 0.0160 | 0.0118 | 0.0290 | 0.0258 |
| 19 | | | 0.7368 | 0.6843 | 0.7313 | 0.0441 | 0.0479 | 0.0349 | 0.1036 | 0.0902 |
| 20 | | | 0.8000 | 0.7507 | 0.7755 | 0.0155 | 0.0295 | 0.0217 | 0.0529 | 0.0471 |
| 21 | | 100 | 0.8000 | 0.7601 | 0.7859 | 0.0129 | 0.0257 | 0.0200 | 0.0529 | 0.0472 |
| 22 | | | 0.8511 | 0.8361 | 0.8452 | 0.0044 | 0.0087 | 0.0071 | 0.0132 | 0.0122 |
| 23 | | | 0.7368 | 0.6825 | 0.7167 | 0.0221 | 0.0453 | 0.0352 | 0.0925 | 0.0827 |
| 24 | | | 0.8000 | 0.7772 | 0.7897 | 0.0077 | 0.0159 | 0.0130 | 0.0244 | 0.0226 |
| 25 | 5 | 30 | 0.2941 | 0.2444 | 0.2614 | 0.0800 | 0.0933 | 0.0854 | 0.0809 | 0.0818 |
| 26 | | | 0.2941 | 0.2551 | 0.2691 | 0.0797 | 0.0899 | 0.0834 | 0.0794 | 0.0770 |
| 27 | | | 0.4000 | 0.3565 | 0.3713 | 0.0654 | 0.0807 | 0.0744 | 0.0702 | 0.0703 |
| 28 | | | 0.4000 | 0.3650 | 0.3768 | 0.0652 | 0.0768 | 0.0722 | 0.0661 | 0.0663 |
| 28 | | 50 | 0.2941 | 0.2697 | 0.2794 | 0.0480 | 0.0529 | 0.0508 | 0.0462 | 0.0465 |
| 30 | | | 0.2941 | 0.2735 | 0.2814 | 0.0478 | 0.0508 | 0.0491 | 0.0456 | 0.0459 |
| 31 | | | 0.4000 | 0.3881 | 0.3962 | 0.0392 | 0.0426 | 0.0409 | 0.0392 | 0.0391 |
| 32 | | | 0.4000 | 0.3800 | 0.3867 | 0.0391 | 0.0431 | 0.0417 | 0.0391 | 0.0391 |
| 33 | | 100 | 0.2941 | 0.2853 | 0.2899 | 0.0240 | 0.0251 | 0.0246 | 0.0235 | 0.0235 |
| 34 | | | 0.2941 | 0.2879 | 0.2917 | 0.0239 | 0.0244 | 0.0240 | 0.0234 | 0.0234 |
| 35 | | | 0.4000 | 0.3897 | 0.3936 | 0.0196 | 0.0206 | 0.0202 | 0.0198 | 0.0197 |
| 36 | | | 0.4000 | 0.3920 | 0.3953 | 0.0196 | 0.0202 | 0.0199 | 0.0196 | 0.0195 |
| 37 | | 30 | 0.7895 | 0.6061 | 0.6614 | 0.0317 | 0.1285 | 0.0588 | 0.6120 | 0.2359 |
| 38 | | | 0.7895 | 0.6474 | 0.6841 | 0.0316 | 0.1070 | 0.0558 | 0.3840 | 0.1761 |
| 38 | | | 0.8000 | 0.6288 | 0.6827 | 0.0287 | 0.1232 | 0.0545 | 0.6135 | 0.2306 |
| 40 | | | 0.8000 | 0.6643 | 0.6986 | 0.0287 | 0.0968 | 0.0515 | 0.3776 | 0.1622 |
| 41 | | 50 | 0.7895 | 0.7491 | 0.7650 | 0.0190 | 0.0363 | 0.0258 | 0.0921 | 0.0501 |
| 42 | | | 0.7895 | 0.7591 | 0.7696 | 0.0190 | 0.0331 | 0.0272 | 0.0434 | 0.0374 |
| 43 | | | 0.8000 | 0.7617 | 0.7763 | 0.0172 | 0.0326 | 0.0239 | 0.0615 | 0.0423 |
| 44 | | | 0.8000 | 0.7706 | 0.7810 | 0.0172 | 0.0288 | 0.0235 | 0.0438 | 0.0343 |
| 45 | | 100 | 0.7895 | 0.7816 | 0.7863 | 0.0095 | 0.0111 | 0.0104 | 0.0117 | 0.0114 |
| 46 | | | 0.7895 | 0.7823 | 0.7861 | 0.0095 | 0.0112 | 0.0105 | 0.0115 | 0.0112 |
| 47 | | | 0.8000 | 0.7920 | 0.7966 | 0.0086 | 0.0104 | 0.0097 | 0.0109 | 0.0106 |
| 48 | | | 0.8000 | 0.7899 | 0.7937 | 0.0086 | 0.0105 | 0.0100 | 0.0106 | 0.0104 |

(1) $\mathcal{S}_3$ is the true value of the consistency in category 3; $\bar{\hat{\mathcal{S}}}_3$ and $\bar{\hat{\mathcal{S}}}_{3U}$ are the sample averages of the 10,000 estimators $\hat{\mathcal{S}}_3$ and $\hat{\mathcal{S}}_{3U}$ respectively.

(2) $V_A\left(\hat{\mathcal{S}}_3\right)$ is the real asymptotic variance, obtained based on the true values of the parameters.

(3) $V_E\left(\hat{\mathcal{S}}_3\right)$ and $V_E\left(\hat{\mathcal{S}}_{3U}\right)$ are the sample variances of the 10,000 estimators $\hat{\mathcal{S}}_3$ and $\hat{\mathcal{S}}_{3U}$ respectively. It is assumed that those sample variances are approximately equal to the true variances of the previous estimators.

(4) $\overline{\hat{V}\left(\hat{\mathcal{S}}_3\right)}$ and $\overline{\hat{V}\left(\hat{\mathcal{S}}_{3U}\right)}$ are sample averages of the 10,000 estimators $\hat{V}\left(\hat{\mathcal{S}}_3\right)$ and $\hat{V}\left(\hat{\mathcal{S}}_{3U}\right)$ respectively.

**Table 6**

**Diagnosis of *n*=100 subjects by two raters in *K*=3 categories (Fleiss *et al.* 2003)**

**(a) Frequencies observed ($x_{ij}$)**

|  | Rater 2 | | | Totals |
|---|---|---|---|---|
| *Rater 1* | *Psychotic* | *Neurotic* | *Organic* | ($x_{i\bullet}$) |
| *Psychotic* | 75 | 1 | 4 | 80 |
| *Neurotic* | 5 | 4 | 1 | 10 |
| *Organic* | 0 | 0 | 10 | 10 |
| *Totals* ($x_{\bullet j}$) | 80 | 5 | 15 | 100 ($x_{\bullet\bullet}$) |

**(b) Parameters estimated under the *delta* model**

| | | | | |
|---|---|---|---|---|
| $\hat{\alpha}_1 = 0.550$ | $\hat{\alpha}_{1U} = 0.575$ | $\hat{\mathcal{S}}_1 = 0.687$ | $\hat{\mathcal{S}}_{1U} = 0.719$ | $\hat{\Delta} = 0.687$ |
| $\hat{\alpha}_2 = 0.037$ | $\hat{\alpha}_{2U} = 0.040$ | $\hat{\mathcal{S}}_2 = 0.500$ | $\hat{\mathcal{S}}_{2U} = 0.528$ | $\hat{\Delta}_U = 0.715$ |
| $\hat{\alpha}_3 = 0.100$ | $\hat{\alpha}_{3U} = 0.100$ | $\hat{\mathcal{S}}_3 = 0.800$ | $\hat{\mathcal{S}}_{3U} = 0.800$ | |

**Table 7**

**Data from a 2×2 table by Nelson and Pepe (2000).**

**(a) Original data (left) and data increased by 0.5 and in a virtual category (right)**

| Categories | 1 | 2 | Total |
|---|---|---|---|
| 1 | 80 | 10 | 90 |
| 2 | 10 | 0 | 10 |
| Total | 90 | 10 | 100 |

| Categories | 1 | 2 | 3 | Total |
|---|---|---|---|---|
| 1 | 80.5 | 10.5 | 0.5 | 91.5 |
| 2 | 10.5 | 0.5 | 0.5 | 11.5 |
| 3 | 0.5 | 0.5 | 0.5 | 1.5 |
| Total | 91.5 | 11.5 | 1.5 | 104.5 |

**(b) Parameters estimated under the *delta* model**

| | | | | |
|---|---|---|---|---|
| $\hat{\alpha}_1^* = 0.680$ | $\hat{\alpha}_{1U}^* = 0.745$ | $\hat{\mathcal{S}}_1 = 0.765$ | $\hat{\mathcal{S}}_{1U} = 0.869$ | $\hat{\varDelta}^* = 0.583$ |
| $\hat{\alpha}_2^* = -0.097$ | $\hat{\alpha}_{2U}^* = -0.031$ | $\hat{\mathcal{S}}_2 = -0.870$ | $\hat{\mathcal{S}}_{2U} = -0.280$ | $\hat{\varDelta}_U^* = 0.714$ |

**(c) Other parameters estimated under the *delta* model if the row rater is a gold standard**

| | |
|---|---|
| $\hat{\mathcal{F}}_1 = \hat{\mathcal{P}}_1 = 0.765$ | $\hat{\mathcal{F}}_{1U} = \hat{\mathcal{P}}_{1U} = 0.839$ |
| $\hat{\mathcal{F}}_2 = \hat{\mathcal{P}}_2 = -0.870$ | $\hat{\mathcal{F}}_{2U} = \hat{\mathcal{P}}_{2U} = -0.280$ |

**Table 8**

**Data from Table IIA by Kramer and Feinstein (1981)**

**(a) Observed frequencies ($x_{ij}$)**

| Rater 1 | Rater 2 | | | | Totals |
|---|---|---|---|---|---|
| | 1 | 2 | 3 | 4 | ($x_{i\bullet}$) |
| 1 | 1 | 2 | 0 | 0 | 3 |
| 2 | 1 | 5 | 3 | 1 | 10 |
| 3 | 1 | 4 | 5 | 2 | 12 |
| 4 | 1 | 1 | 1 | 2 | 5 |
| Totals ($x_{\bullet j}$) | 4 | 12 | 9 | 5 | 30 ($x_{\bullet\bullet}$) |

**(b) Parameters estimated under the *delta* model**

| | | | | |
|---|---|---|---|---|
| $\hat{\alpha}_1 = 0.023$ | $\hat{\alpha}_{1U} = 0.024$ | $\hat{S}_1 = 0.197$ | $\hat{S}_{1U} = 0.206$ | $\hat{\Delta} = 0.182$ |
| $\hat{\alpha}_2 = 0.027$ | $\hat{\alpha}_{2U} = 0.042$ | $\hat{S}_2 = 0.074$ | $\hat{S}_{2U} = 0.115$ | $\hat{\Delta}_U = 0.210$ |
| $\hat{\alpha}_3 = 0.082$ | $\hat{\alpha}_{3U} = 0.092$ | $\hat{S}_3 = 0.234$ | $\hat{S}_{3U} = 0.264$ | |
| $\hat{\alpha}_4 = 0.050$ | $\hat{\alpha}_{4U} = 0.052$ | $\hat{S}_4 = 0.300$ | $\hat{S}_{4U} = 0.311$ | |